\documentclass{m2an}
\usepackage[colorlinks,linkcolor=blue,anchorcolor=blue,citecolor=blue]{hyperref}
\usepackage{placeins}
\usepackage{lipsum}
\usepackage{amsfonts}
\usepackage{graphicx}
\usepackage{epstopdf}
\usepackage{algorithmic}
\usepackage{verbatim}
\usepackage{fancyhdr}
\usepackage{latexsym, epsfig}
\usepackage{pgfplots}
\pgfplotsset{compat=1.18} 
\usepackage{url}
\usepackage{mathrsfs}
\usepackage{multirow}
\usepackage{color}
\usepackage{array}
\usepackage{amsmath}
\numberwithin{equation}{section}
\usepackage{subfigure}
\usepackage{cases}
\usepackage{bm}
\usepackage{float}
\usepackage{ulem}
\usepackage{cancel}

\usepackage[T1]{fontenc}
\usepackage{lmodern}
\usepackage{amssymb}
\usepackage{fix-cm}
\usepackage{tikz}
\usepackage{xparse} 
\usetikzlibrary{shapes.geometric, arrows}
\usetikzlibrary{matrix}
\usepackage[english]{babel} 
\ifpdf
\DeclareGraphicsExtensions{.eps,.pdf,.png,.jpg}
\else
\DeclareGraphicsExtensions{.eps}
\fi

\newcommand{\e}{\varepsilon}

\definecolor{darkred}{rgb}{0.65, 0, 0}

\newcommand{\deletei}{\bgroup\markoverwith{\textcolor{darkred}{\rule[0.5ex]{2pt}{2pt}}}\ULon}
\NewDocumentCommand{\deleteeqi}{ O{darkred} O{2pt} m }{%
  \begingroup
    \sbox0{\ensuremath{#3}}%
    \tikz[baseline=(X.base)]{%
      \node[inner sep=0pt,outer sep=0pt] (X) {\usebox0};%
      \draw[#1,line width=#2] (X.west) -- (X.east);%
    }%
  \endgroup
}

\newcommand{\deleteii}{\bgroup\markoverwith{\textcolor{blue}{\rule[0.5ex]{2pt}{2pt}}}\ULon}
\NewDocumentCommand{\deleteeqii}{ O{blue} O{2pt} m }{%
  \begingroup
    \sbox0{\ensuremath{#3}}%
    \tikz[baseline=(X.base)]{%
      \node[inner sep=0pt,outer sep=0pt] (X) {\usebox0};%
      \draw[#1,line width=#2] (X.west) -- (X.east);%
    }%
  \endgroup
}

\begin{document}
\title{A bi-fidelity method for the uncertain Vlasov-Poisson system near quasineutrality in an asymptotic-preserving particle-in-cell framework} \thanks{Submitted to the editors January 16, 2025.}

\author{Guangwei Liu}\address{Beijing Computational Science Research Center, Beijing, 100193, China;
	\email{guangweiliu@csrc.ac.cn\ \&\ ylwang@csrc.ac.cn}}
\author{Liu Liu}\address{Department of Mathematics, The Chinese University of Hong Kong, Hong Kong;
\email{lliu@math.cuhk.edu.hk}}
\author{Yanli Wang}\sameaddress{1}
%
%

\begin{abstract} 
	   In this paper, we study the Vlasov-Poisson system with massless electrons (VPME) near quasineutrality and with uncertainties. Based on the idea of reformulation on the Poisson equation by [P. Degond et.al., \textit{Journal of Computational Physics},  229 (16), 2010, pp. 5630--5652], we first consider the deterministic problem and develop an efficient asymptotic preserving particle-in-cell (AP-PIC) method to capture the quasineutral limit numerically, without resolving the discretizations subject to the small Debye length in plasma. The main challenge and difference compared to previous related works is that we consider the nonlinear Poisson in the VPME system, which contains $e^{\phi}$ (with $\phi$ being the electric potential) and provide an explicit scheme. In the second part, we extend to study the uncertainty quantification (UQ) problem and develop an efficient bi-fidelity method for solving the VPME system with multidimensional random parameters, by choosing the Euler-Poisson equation as the low-fidelity model. Several numerical experiments are shown to demonstrate the asymptotic-preserving property of our deterministic solver and the effectiveness of our bi-fidelity method for solving the model with random uncertainties.   
\end{abstract}
\maketitle
\section{Introduction}	
	Kinetic equations of Vlasov type are in widespread use as models in plasma physics. 
	In plasma simulations, the basic kinetic model is the Vlasov equation, coupled with electromagnetic field equations. In this work we consider the model for collisionless, unmagnetised plasma. In particular, we will study the Vlasov-Poisson system with massless electrons, namely the VPME system: 
	\begin{equation}
		\label{Intro:VP-semi}
		\text{(VPME)}: 
		\left\{\begin{array}{l}
			\displaystyle 
			\partial_t f + \mathbf{v}\cdot \nabla_{\mathbf{x}} f  - \nabla_{\mathbf{x}} \phi\cdot \nabla_{\mathbf{v}} f =0, \\[6pt]
			-\e^2 \Delta_{\mathbf{x}} \phi  = \int f \mathrm{d}\mathbf{v} - e^{\phi}, \\[6pt]
			f(0,\mathbf{x},\mathbf{v}) = f_0, \, \int f_0 \, \mathrm{d}\mathbf{v}\mathrm{d}\mathbf{x} = 1.  
		\end{array}\right.
	\end{equation}
	Here $\e$ is the dimensionless parameter characterized by the Debye length and $f(t, \mathbf{x},\mathbf{v})$ represents the distribution of ions. This VPME system has been used in plasma physics to model ion plasma, study the formation of ion-acoustic shocks \cite{SCM} and the expansion of plasma into vacuum \cite{Med}, among other applications. We refer to \cite{GP} for more detailed introduction.

	In real applications of plasma physics, these equations often exhibit different scales, ranging from the microscopic dynamics of individual particles to the macroscopic behavior of the plasma as a whole. When the scaled Debye length $\e\ll 1$, the plasma is called {\it quasineutral}: since the scale of charge separation is very small, the plasma appears to be neutral at the scale of observation. Quasineutrality is a common and important phenomenon for plasmas. In \cite{Chen} the author emphasizes quasineutrality as one of the key properties that distinguish plasmas from ionised gases, we also refer to \cite{Review} and other theoretical works \cite{bouchut1991global, Brenier2000, BG1, crouseilles2016multiscale, Golse, G1, G2, griffin2020global, griffin2021global, han2011quasineutral,  han2017quasineutral} which give a review of recent developments on for the Vlasov-type equations near quasineutrality. 
	
	For numerical computations, solving these equations with small physical parameter can be significantly expensive and challenging. Asymptotic-preserving (AP) schemes, originated from the work of S. Jin for multiscale kinetic problems \cite{jin1999efficient, jin1991discrete, jin1993fully}, have been popularly used to tackle this difficulty. Denote $\mathcal{F}_{\varepsilon}$ as the Vlasov-Poisson system and $ \mathcal{F}_{0}$ as the corresponding quasineutral Vlasov equation, i.e., the limiting system of $\mathcal{F}_{\e}$ as $\e\to 0$. The AP scheme requires that, without resolving the space and time steps associated with $\e$, when $\e\to 0$, it automatically becomes a consistent and stable discretization of $\mathcal{F}_0$. In this work, one of our goals is to develop an AP scheme for the VPME system \eqref{Intro:VP-semi} that can capture the quasineutral limit on the discretized level. 
	
	Another computational challenges of the Vlasov-type equations arise from the high-dimensionality. Particle methods are preferred in the simulation of plasma kinetic models--given as 7-dimensional problem (3D in space, 3D in velocity, plus time), in comparison to the grid-based Eulerian methods. In Particle-in-Cell (PIC) methods, the coupling between the particles and the field is implemented through a space grid. Upon solving the positions and velocities along the characteristic lines, one can adopt finite-difference methods to obtain the field solution on the grid, then interpolate back to the positions of particles. We mention two books \cite{birdsall2018plasma,hockney2021computer} for an overview on this topic. The convergence of PIC methods has also been investigated \cite{Ganguly1989}. 
	
	Numerical computation for the VPME system \eqref{Intro:VP-semi}, that is under study in this work, seems underdeveloped. For the linear Poisson equation, a key idea was first proposed by P. Degond et al. \cite{crispel2005asymptotically}, which is to reformulate the system to a new one that remains uniformly elliptic thus can capture the quasineutral limit of the system. Later in \cite{degond2010asymptotic}, the authors developed an asymptotic-preserving PIC method for the Vlasov-Poisson system that can capture the quasineutral limit at the discrete level. Previous works on AP methods for the quasineutral limit have been devoted to the Euler-Poisson model \cite{crispel2007, crouseilles2024high, Degond2008} and the Vlasov-Poisson problem \cite{Crouseilles2009}. In this work, we will develop an efficient AP-PIC method for solving the VPME system \eqref{Intro:VP-semi}. For the deterministic problem, the main difficulty compared to the previous work \cite{degond2010asymptotic} is that we consider the nonlinear Poisson equation with $e^{\phi}$ on the right-hand-side, thus the derivation of the quasineutral limit and the design of our discretized scheme is {\it novel and new}. Moreover, we provide a stable, explicit scheme that allows us to solve the reformulated Poisson equation efficiently. 
	
	In industrial applications of plasma, it is important to study the uncertainty quantification (UQ) problems since there are many sources of uncertainties in the model. In the second part of this paper, we will extend to study UQ problems for the VPME system and consider the non-intrusive stochastic collocation (SC) method \cite{Babuska, Schwab, Webster,Xiu2005} in a bi-fidelity framework. One key challenge of SC approaches is the computational cost, as they require repetitive implementations of the deterministic solver. Fortunately, in practice we can usually find some approximated low-fidelity models that are computationally cheaper but can still characterize behaviour of the complex model in the random space. The multi-fidelity algorithm \cite{Narayan2014,ZNX2014} is to combine the computational efficiency of low-fidelity models with the high accuracy of high-fidelity models, to construct an accurate surrogate for the high-fidelity model at a significantly reduced simulation cost. The feasibility of multi-fidelity methods and empirical error estimations to predict whether the multi-fidelity approach works or not for practical problems, has been discussed in \cite{Gao-Zhu-Wang-2020}. Regarding UQ for kinetic equations with random parameters, one can refer to a series of works using the bi-fidelity methods \cite{BLPZ2022,LPZ2022,LZ2020,pareschi2021introduction,LinLiu2025}. 

	The novelty and main contribution of our work is two-fold: (i) In this work, inspired by \cite{degond2010asymptotic}, where the reformulation idea is introduced, we design a PIC method for the VPME system that includes a semi-implicit particle update and an implicit reformulated Poisson equation. Our proposed scheme is asymptotic-preserving and can capture the correct quasi-neutral limit, while allowing the time step and mesh size to be much larger than the Debye length. Moreover, our AP-PIC scheme, built based on reformulation, decouples the field-particle iteration, where the particle trajectories and macroscopic moments are computed only once per time step, significantly improving computational efficiency compared with a fully implicit PIC scheme \cite{chen2011energy}. We also mention that compared with the work by Degond et al. \cite{P2012Numerical} that considers the macroscopic system, our model involves the Vlasov equation on the kinetic level, thus their implicit scheme can not be directly adapted here. (ii) To our best knowledge, this is the first work to study the UQ problem for the VPME system. We explore different multi-fidelity approaches using appropriately chosen low-fidelity models and evaluate their accuracy and efficiency. For the high-fidelity solver, the AP-PIC method developed in (i) is adopted, which further underscores the importance of achieving high efficiency for our deterministic solver.
   
   The rest of the paper is organized as follows. In Section \ref{sec:models}, we first give an introduction of the Vlasov-Poisson system with massless electron and its related macroscopic equations. In Section \ref{sec:AP-PIC}, we begin with reviewing the classical PIC method, then derive the reformulated Poisson equation that leads to the design of our fully discretized AP-PIC scheme, with a summary of the method shown in Section \ref{fig:flowchart}. In Section \ref{sec:UQ}, we introduce the bi-fidelity method and how we choose the high-fidelity and low-fidelity model (Euler-Poisson) in solving the VPME system with uncertainties. Several numerical examples, including the deterministic and problem with random parameters have been shown in Section \ref{sec:Num}, which demonstrate the effectiveness and AP property of our solver for the VPME system, in addition to accuracy and efficiency of our designed bi-fidelity method for the corresponding UQ problems. Finally, we conclude the paper and mention some future work. 
\section{Models}
\label{sec:models}

We are interested in the systems of Vlasov-Poisson type, which are kinetic equations describing dilute, collisionless, weakly magnetised plasmas. Combining the classical Vlasov-Poisson system with a Maxwell-Boltzmann law for the electron distribution leads to the Vlasov-Poisson system with massless electrons \cite{Review,han2011quasineutral}, or the VPME system, with the dimensionless form given by 
\begin{subnumcases}{\label{VP-semi} \text{(VPME):}}
	\partial_t f + \mathbf{v}\cdot \nabla_{\mathbf{x}} f - \nabla_{\mathbf{x}} \phi\cdot \nabla_{\mathbf{v}} f =0, \label{eq:VPME-Vlasov} \\[2pt] 
	-\e^2 \Delta_{\mathbf{x}} \phi  = \int f \mathrm{d}\mathbf{v} - e^{\phi},\label{eq:VPME-Poisson}\\[2pt]
	f(0,\mathbf{x},\mathbf{v}) = f_0, \, \int f_0 \, \mathrm{d}\mathbf{v}\mathrm{d}\mathbf{x} = 1.  \label{eq:VPME-init}
\end{subnumcases}
where $f(t,\mathbf{x},\mathbf{v})$ is interpreted as the distribution of ions, at position $\mathbf{x}$ moving with velocity $\mathbf{v}$. 
Here $\phi$ is the electric potential that is governed by the Poisson equation, and $\e$ is the dimensionless parameter characterized by the Debye length. More details of derivation of the model can be found in \cite{han2011quasineutral}. 

Based on the relation \eqref{VP-semi} and we assume the solution $f$ given by the monokinetic form: 
\begin{equation}
	\label{eq:mono} 
	f(t,\mathbf{x},\mathbf{v}) = n(t,\mathbf{x})\delta(\mathbf{v}-\mathbf{u}(t,\mathbf{x})), 
\end{equation}
with $\delta$ being the Dirac delta function. Here, $n(t,\mathbf{x})$ is the ion density and $u(t,\mathbf{x})$ is the ion velocity given by
\begin{equation}
	\label{ion-density-velocity}
	n(t,\mathbf{x}) = \int f(t,\mathbf{x},\mathbf{v})\mathrm{d}\mathbf{v},\quad \mathbf{u}(t,\mathbf{x}) = \frac{1}{n(t,\mathbf{x})}\int \mathbf{v}f(t,\mathbf{x},\mathbf{v})\mathrm{d}\mathbf{v}.
\end{equation}
Then one derives {\it isothermal compressible Euler system} (ICE) by taking the first and second moments in velocity: 
\begin{equation}
	\label{Euler-semi}
	\text{(ICE): } \left\{\begin{array}{l}
		\partial_{t} n+\nabla_{\mathbf{x}} \cdot(n \mathbf{u})=0 ,\\[6pt]
		\partial_{t} (n\mathbf{u})+ \nabla_{\mathbf{x}} \cdot (n\mathbf{u} \otimes \mathbf{u})=- \nabla_{\mathbf{x}} n. 
	\end{array}\right.
\end{equation}

For the VPME system \eqref{VP-semi}, periodic boundary conditions are assumed for the distribution $f(t,\mathbf{x}, \mathbf{v})$ and potential $\phi$. 
In the one-dimensional case, the boundary conditions are given by
\begin{equation}
	\label{model_bc}
	\begin{aligned}
		&  f(t, x_L, v) = f(t, x_R ,v), \qquad
		\phi(t, x_L) = \phi(t, x_R). \\[4pt]
	\end{aligned}
\end{equation}

\section{AP-PIC method for VPME system}
\label{sec:AP-PIC}
\subsection{Classical PIC method}
\label{subsec:classical-pic}
The numerical simulation of the Vlasov-type equations is usually performed by particle-in-cell (PIC) methods, which approximate the plasma by a relatively small number of particles. The trajectories of these particles are computed from the characteristic curves corresponding to the Vlasov equation, and the self-consistent electric field is computed on a mesh of the physical space. We refer to some textbooks \cite{birdsall2018plasma,grigoryev2012numerical,hockney2021computer} for an exposition of PIC methods. Below we consider one-dimensional spatial and velocity variables.
In the standard particle method, one approximates the initial distribution $f_0$ by the following 
Dirac mass sum: 
$$ f_{N_{\text{total}}}^0(x,v) = \sum_{k=1}^{N_{\text{total}}} w_k \delta(x-X_k) \delta(v-V_k), $$
where $(X_k^0, V_k^0)_{1\leq k\leq N_{\text{total}}}$ is a beam of $N_{\text{total}}$ particles distributed in the phase space, 
$w_{k}$ represents the weight of the $k$-th particle.
Here we divide the space $\left[x_{L}, x_{R}\right]$ with a uniform space step $\Delta x = h = \frac{x_{R}-x_{L}}{N_{h}}$ and $x_{j} = (j-\frac{1}{2})h$, for $j = 1,\dots,N_{h}$ represent the mesh center point. 

We refer to Appendix for some details on the initial charge assignment. In numerical simulations, the Dirac mass of velocity is usually replaced by a smooth function, which is chosen by the  
cosine kernel function mentioned in \cite{engquist2003dirac, jin2005computing}. In particular, we let

\begin{subequations}
	\begin{align}
		& f_{{N_{\text{total}}},\alpha}(0,x,v) := \sum_{k=1}^{N_{\text{total}}} w_k W(x-X_k(0)) B_{\eta}(v-V_k(0)),  \label{eq:f-PIC} \\
		& B_{\eta}(v-V_k(0)) =\left\{\begin{array}{ll}
			\frac{1}{2 \eta}\left(1+\cos \frac{|\pi \left(v-V_k(0)\right)|}{\eta}\right), & \Big|\frac{v-V_k(0)}{\eta}\Big| \leqslant 1, \\[8pt]
			0, & \Big|\frac{v-V_k(0)}{\eta}\Big|>1, 
		\end{array}\right.  \label{eq:delta-v-approx}
	\end{align}
\end{subequations}
where $\eta$ is the regularization parameter, and 
\begin{equation}
	\label{eq:pic-assignment-function}
	W(x)=\left\{\begin{array}{ll}
		1 & |x|<h / 2 \text { or } x=h / 2, \\[4pt]
		0 & \text { otherwise. }
	\end{array}\right.
\end{equation}
As time evolves, one approximates the solution $f$ by 
$$ f_{N_{\text{total}}}(t,x,v) = \sum_{k=1}^{N_{\text{total}}} w_k W(x-X_k(t)) B_{\eta}(v-V_k(t)), $$
where $(X_k, V_k)_{1\leq k\leq N_{\text{total}}}$ is the position in phase space of particle $k$ moving along the characteristic curve defined by 
\begin{equation}
	\left\{
	\begin{array}{ll}
		\frac{dX}{dt} = V,  \\[6pt]
		\frac{dV}{dt} = -\nabla\phi,  \\[6pt]
		X(t=0)=x^0, \quad V(t=0)=v^0 . 
	\end{array}
	\right.
\end{equation}
The volume average of density, current density and momentum flux is approximated by
\begin{subequations}
	\begin{align}
		& n(x,t) = \frac{1}{h}\int_{x_{j}-\frac{h}{2}}^{x_j+\frac{h}{2}}\int_{-\infty}^{+\infty}f(t,x,v)\mathrm{d}v \mathrm{d}x\approx \frac{1}{h}\sum_{k=1}^{N_{\text{total}}}w_{k} W(x-X_{k}(t)), \label{eq:moment0-approx} \\
		& J(x,t) = \frac{1}{h}\int_{x_{j}-\frac{h}{2}}^{x_j+\frac{h}{2}}\int_{-\infty}^{+\infty}v f(t,x,v)\mathrm{d}v\mathrm{d}x \approx \frac{1}{h}\sum_{k=1}^{N_{\text{total}}}w_{k}V_{k}(t)W(x-X_{k}(t)), \label{eq:moment1-approx} \\
		& S(x,t) = \frac{1}{h}\int_{x_{j}-\frac{h}{2}}^{x_j+\frac{h}{2}}\int_{-\infty}^{+\infty}v^{2}f(t,x,v)\mathrm{d}v\mathrm{d}x \approx \frac{1}{h}\sum_{k=1}^{N_{\text{total}}}w_{k}\left(V_{k}(t)\right)^{2}W(x-X_{k}(t)). \label{eq:moment2-approx}
	\end{align}
\end{subequations}
where $J(x,t) = \left(nu\right)(x,t)$ and $x\in \left[x_j - \frac{h}{2}, x_j+\frac{h}{2}\right],\: j = 1,\dots,N_{h}$.

In the classical PIC method, given the information of particles at time $t^m = m\Delta t$, we solve the following ODE system and update the positions and velocities at time $t^{m+1}$: 
\begin{equation}
	\label{Char}
	\left\{
	\begin{array}{ll}
		\frac{dX_k}{dt} = V_k,  \\[4pt]
		\frac{dV_k}{dt} = -\nabla\phi(t,X_k),  \\[4pt]
		X_k(t^m)=x_k^m, \quad V_k(t^m)=v_k^m. 
	\end{array}
	\right.
\end{equation}
In the remaining of this section, we will study how to compute $\phi$ by the reformulated Poisson equation and introduce a fully discretized scheme that satisfies the asymptotic-preserving property for the VPME system \eqref{VP-semi}. 
\subsection{Derivation of the reformulated Poisson equation}
\label{subsec:RF-Poisson}
In this section, we will introduce the derivation of the reformulated Poisson equation in order to design numerical scheme to capture the quasineutral limit of the VPME system. 

For the Vlasov-Poisson system in a nearly quasineutral regime, the classical PIC method suffers from stability constraints on the time and space steps related to the small Debye length and large plasma frequency. Inspired by the work in \cite{degond2010asymptotic} where a reformulation of the Poisson equation (constituting both electrons and ions) is introduced, we propose an Asymptotic-Preserving PIC method for the VPME system \eqref{VP-semi}, which provides a consistent PIC discretization of the quasineutral Vlasov equation.

We begin from the mass and momentum conservation equations:
\begin{equation}
	\label{M1}
	\left\{\begin{array}{l}
		\partial_{t} n_{\alpha}+\nabla \cdot(n u)_{\alpha}=0, \\[4pt]
		\partial_{t}(n u)_{\alpha}+\nabla \cdot S_{\alpha}=-n_{\alpha} \nabla \phi, 
	\end{array}\right. 
\end{equation}
where $\alpha=i,e$, $(nu)_{\alpha}$ are the ion and electron momenta and $S_{\alpha}$ are the momentum fluxes
\begin{equation}\label{moments}
	(n u)_{i, e}=\int_{\mathbb{R}^{d}} f_{i, e}(x, v, t) v \mathrm{d} v, \quad S_{i, e}=\int_{\mathbb{R}^{d}} f_{i, e}(x, v, t) v \otimes v \mathrm{d} v. 
\end{equation}
and the symbol $\otimes$ denotes the tensor product. \\
Take the time derivative of the mass conservation equation for ion,  we have
\begin{equation}\label{ion_mass_conv}
	\partial_{tt} n_{i} +\nabla \cdot \partial_{t} (n u)_{i} =0. 
\end{equation}
Replacing the term $\partial_{t} (n u)_{i}$ by using the momentum equation in \eqref{M1} for $\alpha=i$, one gets
\begin{equation}\label{ni_eqn}
	\partial_{tt} n_{i} - \nabla^{2}:S_{i} = \nabla \cdot (n_{i}\nabla \phi).
\end{equation}
Take second order derivative in time on the Poisson equation, 
\begin{equation}\label{P-time}
	\partial_{tt}\left(-\varepsilon^2 \Delta\phi  - (n_{i} - n_{e})\right) = 0.
\end{equation}
substituting the term 
$\partial_{tt}n_i$ in \eqref{ni_eqn} by using \eqref{P-time}, one derives
\begin{equation}\label{P2}
	-\nabla \cdot\left( (\varepsilon^2 \partial_{tt} + n_{i} ) \nabla \phi \right) + \partial_{tt} n_{e} = \nabla^{2}:S_{i}, 
\end{equation}
where $\nabla^2$ denotes the tensor of second order derivatives and “:” is the contracted product of two tensors. 

If the initial data $\phi_{0} := \phi|_{t=0}$ and $\phi^{\prime}_{0}:=\partial_{t}\phi|_{t=0}$ satisfy the following Poisson equations at the initial time:
\begin{equation}\label{IC-phi}
	\begin{array}{l}
		-\varepsilon^{2} \Delta \phi_{0}=\left(n_{i}-n_{e}\right)_{0} , \\[4pt]
		-\varepsilon^{2} \Delta \phi_{0}^{\prime}=\left(n_{i}-n_{e}\right)_{0}^{\prime}=-\nabla\cdot\left((nu)_{i} - (nu)_{e}\right)_{0} , 
	\end{array}
\end{equation}
with $n_e=e^{\phi}$, the original Poisson equation in our VPME system
\begin{equation}\label{Poisson}
	-\e^2 \Delta \phi = n_i - e^{\phi}, 
\end{equation}
is satisfied at all time. 
By replacing $n_e$ by $e^{\phi}$ in \eqref{P2}, we summarize the reformulated Poisson equation as
\begin{equation}\label{RF-Poisson}
	-\nabla \cdot\left((\varepsilon^2 \partial_{tt} + n_{i} ) \nabla \phi \right) + \partial_{tt} e^{\phi} = \nabla^{2}:S_{i}.
\end{equation}
Thus the Vlasov-Poisson system with massless electrons (VPME) \eqref{VP-semi} is equivalent to the reformulated VPME system as follows,
\begin{eqnarray}
	\partial_{t}f + v \cdot \nabla_{x} f -\nabla_{x} \phi \cdot \nabla_{v} f=0, \label{RF_VPME1}\\[4pt]
	-\nabla \cdot\left((\varepsilon^2 \partial_{tt} + n_{i} ) \nabla \phi \right) + \partial_{tt} e^{\phi} = \nabla^{2}:S_{i}.\label{RF_VPME2}
\end{eqnarray}
where $\phi$ need to satisfy the initial data \eqref{IC-phi}.
As $\e\to 0$, the quasineutral Vlasov system becomes 
\begin{eqnarray}
	\partial_{t}f + v \cdot \nabla_{x} f -\nabla_{x} \phi \cdot \nabla_{v} f=0, \label{Quasi-Limit1} \\[4pt]
	-\nabla \cdot\left( n_i \nabla\phi \right) + \partial_{tt} (e^{\phi}) = \nabla^{2}:S_{i}, \label{Quasi-Limit2}
\end{eqnarray}
together with the initial conditions \eqref{IC-phi}.

We note that although the original Poisson equation \eqref{Poisson} and the reformulated Poisson equation \eqref{RF-Poisson} are equivalent, they are not equally suited in the quasi-neutral limit.

\subsection{The AP-PIC method}
\label{subsec:scheme}
First, we design a semi-implicit time-discretization scheme for the reformulated Poisson equation \eqref{RF-Poisson}: 
\begin{equation}\label{semi-scheme}
	\begin{array}{l}
		- \nabla \cdot \Big( (\Delta t)^2 n_i^m \nabla\phi^{m+1} + 
		\e^2 \left(\nabla\phi^{m+1} - 2 \nabla\phi^{m} + \nabla\phi^{m-1}\right) \Big)  \\[4pt]
		= (\Delta t)^2 \nabla^{2}:S_{i}^m - ( e^{\phi^{m+1}} - 2 e^{\phi^{m}} + e^{\phi^{m-1}}). 
	\end{array}
\end{equation}
Note that for both the term $\nabla \cdot ( \partial_{tt}\nabla\phi )$ and $\partial_{tt}e^{\phi}$ we use an implicit disretization, since it would be convenient to replace the term $e^{\phi^{m}}$ and $e^{\phi^{m-1}}$ through the \eqref{Poisson}, which means
\begin{equation}
	\label{RF-Poisson-time-discretization}
	- \nabla \cdot \Big( (\Delta t)^2 n_i^m \nabla\phi^{m+1}+\e^2\nabla\phi^{m+1}\Big) +  e^{\phi^{m+1}} \\[4pt]
	= (\Delta t)^2 \nabla^{2}:S_{i}^m + 2 n_i^m - n_i^{m-1}. 
\end{equation}
In the limit $\e\to 0$, one has
\begin{equation}
	\label{Quasi-Limit-time}
	- \nabla \cdot \Big( (\Delta t)^2 n_i^m \nabla\phi^{m+1} \Big) + e^{\phi^{m+1}}
	= (\Delta t)^2 \nabla^{2}:S_{i}^m + 2n_i^{m} - n_i^{m-1}.
\end{equation}
which is a consistent discretization of the quasineutral potential equation \eqref{Quasi-Limit2}. Thus, our scheme \eqref{semi-scheme} satisfies the property of {\it Asymptotic-Preserving} (AP). 

Consider the right-hand side of \eqref{RF-Poisson-time-discretization}, the term $n_i^m-n_i^{m-1}$ can be replaced by $-\nabla(nu)_i^m$, we derive an equivalent form shown as
\begin{equation}
	\label{RF-Poisson-time-discretation-2}
	- \nabla \cdot \Big( (\Delta t)^2 n_i^m \nabla\phi^{m+1}+\e^2\nabla\phi^{m+1}\Big) +  e^{\phi^{m+1}} \\[4pt]
	= (\Delta t)^2 \nabla^{2}:S_{i}^m + n_i^m - \nabla (nu)_i^{m}. 
\end{equation}
Let $\Delta x = h = \frac{x_{R}-x_{L}}{N_{h}}$ the uniform space step and $x_j=(j-\frac{1}{2})h$, for $j=1,\cdots, N_{h}$ the center of each cell. 
Denote $g_j\approx g(x_j)$ as the approximated volume average of function $g$ on the space grid, with $x_{j} \in [x_{j-\frac{h}{2}}, x_{j+\frac{h}{2}}],\: j=1,2,\dots,N_{h}$. For the space discretizations, we set 
$$
\begin{array}{l}
	\left(D^{+} g\right)_{j}=\frac{g_{j+1}-g_{j}}{h}, \\[4pt]
	\left(D^{-} g\right)_{j}=\frac{g_{j}-g_{j-1}}{h}, \\[4pt]
	\left(\Delta g\right)_{j}=\left(D^{+} D^{-} g\right)_{j}=\left(D^{-} D^{+} g\right)_{j}
	=\frac{g_{j+1}-2 g_{j}+g_{j-1}}{h^{2}}.
\end{array}
$$
A full discretization of the reformulated Poisson equation \eqref{RF-Poisson-time-discretation-2} in the 1D case is given by 
\begin{equation}
	\label{Full-scheme}
	\begin{array}{l}
		-\e^2 \left(\Delta\phi^{m+1}\right)_{j}- (\Delta t)^2 \Big(D^{-}\left(n_{i}^{m}D^{+} \phi^{m+1}\right)\Big)_{j} + e^{\phi_j^{m+1}}  \\[4pt]
		= (\Delta t)^2 \left(\Delta S_{i}^{m}\right)_{j}
		+  (n_i)_{j}^m -D^{-}(n_{i}u_{i})_j^m. 
	\end{array}
\end{equation}
Due to the nonlinear term $e^{\phi_j^{m+1}}$, the Newton iteration method is used to solve this scheme. 

Define the coefficient matrix $A^{m}$ by 
$$
\begin{aligned}
	\left(A^{m} \phi\right)_{j}= & -\e^2\left(\Delta \phi\right)_{j}
	- (\Delta t)^2 \Big(D^{-}\left(n_{i}^{m} D^{+}\phi\right)\Big)_{j} \\[4pt]
	= & \left(\frac{2 \e^2}{h^{2}} + \frac{(\Delta t)^2\left(n_{i}^{m}\right)_{j}}{h^{2}}+\frac{(\Delta t)^2\left(n_{i}^{m}\right)_{j-1}}{h^{2}}\right)\phi _{j} \\[2pt]
	& -\left(\frac{\varepsilon^{2}}{h^{2}}+\frac{(\Delta t)^2\left(n_{i}^{m}\right)_{j}}{h^{2}}\right) \phi_{j+1} \\[2pt]
	& -\left(\frac{\varepsilon^{2}}{h^{2}}+\frac{(\Delta t)^2\left(n_{i}^{m}\right)_{j-1}}{h^{2}}\right) \phi_{j-1},
\end{aligned}
$$
then \eqref{Full-scheme} is transformed to solving the following problem, 
\begin{subequations}
	\begin{equation}
		\label{F_phi_b}
		F(\mathbf{\Phi}^{m+1})= A^{m}\mathbf{\Phi}^{m+1} + \mathrm{diag}(e^{\mathbf{\Phi}^{m+1}}) - \mathbf{b}^m = 0,
	\end{equation}
	\begin{equation}
		\label{Phi_vector}
		\mathbf{\Phi}^{m+1} = \left[\phi^{m+1}_{1},\phi^{m+1}_{2},\dots,\phi^{m+1}_{N_{h}}\right]^\top, 
	\end{equation}
	\begin{equation}
		\label{Full-scheme-right}
		b_{j}^m = (\Delta t)^2 \left(\Delta S_{i}^{m}\right)_{j}
		+  (n_i)_{j}^m - \left( (n_iu_i)_{j}^{m} - (n_iu_i)_{j-1}^{m} \right), \qquad \text{for } j=1,2,\cdots,N_{h}. 
	\end{equation}
\end{subequations}
We consider periodic boundary condition \eqref{model_bc} for $\phi$, thus
\begin{equation}
	\label{appic_bc}
	\phi_{0} = \phi_{N_{h}},\quad \phi_{1} = \phi_{N_{h}+1}. 
\end{equation}
Then the specific form of $A^m$ is derived as
\begin{equation}
	\label{appic_matrix}
	A^{m} =
	\begin{bmatrix}
		a_{0} + a_{1} & -a_{1} & 0 & \cdots & -a_{0} \\
		-a_{1} & a_{1} + a_{2} & -a_{2} & \cdots & 0 \\
		0 & -a_{2} & a_{2}+a_{3} & \cdots & 0 \\
		\vdots & \vdots & \vdots & \ddots & \vdots \\
		0 & 0 & -a_{N_{h}-2} & a_{N_{h}-2} + a_{N_{h}-1} & -a_{N_{h}-1}\\
		-a_{N_{h}} & 0 & 0 & -a_{N_{h}-1}& a_{N_{h}-1} + a_{N_{h}}
	\end{bmatrix}, 
\end{equation}
where
$$
a_{j} = -\frac{\varepsilon^{2}}{h^{2}} - \frac{(\Delta t)^{2} \left(n_{i}^{m}\right)_{j}}{h^{2}}, \qquad k= 1,2,\dots,N_{h}.
$$
Here $a_{0}$ is set based on the periodic boundary condition for $f$, namely 
\begin{eqnarray}
	\label{eq:bc-f}
	a_{0} = a_{N_{h}} = -\frac{\varepsilon^{2}}{h^{2}} - \frac{(\Delta t)^{2} \left(n_{i}^{m}\right)_{N_{h}}}{h^{2}}.
\end{eqnarray}
The Newton iteration needs the Jacobian of $F(\mathbf{\Phi}^{m+1})$ defined in \eqref{F_phi_b}, 
which is
\begin{equation}
	\label{Jacobian_F}
	J(\mathbf{\Phi}^{m+1}) = A^{m} + \mathrm{diag}(e^{\mathbf{\Phi}^{m+1}}) . 
\end{equation}
Then the solution is updated by
\begin{equation}
	\label{Newton-method-update}
	\mathbf{\Phi}^{m+1,l+1} = \mathbf{\Phi}^{m+1,l}- J(\mathbf{\Phi}^{m+1,l})^{-1}F(\mathbf{\Phi}^{m+1,l}),
\end{equation}
where $l$ represents the number of iterations. We set the termination condition as
\begin{equation}
	\label{eq:Newton-iter-condition}
	\left \| \mathbf{\Phi}^{m+1,l} - \mathbf{\Phi}^{m+1,l} \right \|_{L^{2}} \le \eta_{1}. 
\end{equation}
where $\eta_{1}$ is the small error tolerance, and the norm is $L^{2}$ norm.  
Let $\mathbf{\Phi}^{m+1,\ast}$ be the convergent solution of the Newton method, then we update $\mathbf{\Phi}^{m+1}= \mathbf{\Phi}^{m+1,\ast}$. Finally, the electric field $ -\nabla \mathbf{\Phi}^{m+1}$ 
can be approximated using the first-order central difference scheme, i.e.,
\begin{equation}
	\label{eq:electric_field}
	- \nabla \phi_{j} \approx -\frac{\phi_{j+1} - \phi_{j-1}}{2h},\qquad j = 1,2,\dots,N_{h}.
\end{equation}

\bigskip

\noindent\textbf{The discretized AP-PIC scheme.}
We now introduce the time discretization for the particle positions and velocities shown in \eqref{Char} and summarize our AP-PIC method for solving the VPME system \eqref{RF_VPME1}--\eqref{RF_VPME2}. Let the time-stepping method for the position equation be fully implicit, while that of the velocity equation be semi-implicit: the electric field is computed implicitly at the particle positions obtained from previous time step. 
Assume at time $t^m=m\Delta t$, the particles are located in $(x_k^m, v_k^m)$ for $1\leq k\leq N_{\text{total}}$,  
then the time advancement of our scheme is given by 
\begin{equation}
	\label{partcle_discretization}
	\left\{\begin{array}{l}
		\frac{V_k^{m+1} - V_k^m}{\Delta t} = -\nabla\phi^{m+1}(X_k^m), \\[6pt]
		\frac{X_k^{m+1} - X_k^m}{\Delta t} = V_k^{m+1}, \\[6pt]
		X_k^m = x_k^m, \qquad V_k^m = v_k^m, 
	\end{array}\right. 
\end{equation}
We summarize our proposed AP-PIC method as a flowchart in the following Figure \ref{fig:flowchart}, which shows how to update the variables from time \( t^{m} \) to \( t^{m+1} \).
\tikzstyle{startstop} = [rectangle,rounded corners, minimum width=3cm,minimum height=1cm,text centered,text width =6cm, draw=black,fill=red!30]
\tikzstyle{io} = [trapezium, trapezium left angle = 70,trapezium right angle=110,minimum width=6cm,minimum height=1cm,text centered,text width =6cm,draw=black,fill=blue!30]
\tikzstyle{process} = [rectangle,minimum width=6cm,minimum height=1cm,text centered,text width =8cm,draw=black,fill=orange!30]
\tikzstyle{decision} = [diamond,aspect = 3,text centered,draw=black,fill=green!30]
\tikzstyle{arrow} = [thick,->,>=stealth]
\tikzstyle{straightline} = [line width = 1pt,-]
\tikzstyle{point}=[coordinate]
\begin{figure}[!hptb]
	\centering
	\begin{tikzpicture}[node distance=2cm]
		\node (input1) [io] {
			Velocities, positions of particles $\{\left(X_{k}^{m}, V_{k}^{m}\right)\}_{k=1}^{N_{\text{total}}}$ and \\ macroscopic quantities $\mathbf{n}^{m}$, $\mathbf{J}^{m}$, $\mathbf{S}^{m}$ }; 
		\node (process1) [process,below of=input1] {Compute $\mathbf{\Phi}^{m+1}$ by the reformulated Poisson equation using \eqref{Full-scheme} with Newton's method};
		\node (process2) [process,below of=process1] {Compute electric field $-\nabla \mathbf{\Phi}^{m+1}$ through \eqref{eq:electric_field}};
		\node (process3) [process,below of=process2] {Update velocities and positions of particles $\{\left(X_{k}^{m+1}, V_{k}^{m+1}\right)\}_{k=1}^{N_{\text{total}}}$ by the ODE system \eqref{partcle_discretization}};
		\node (stop) [startstop,below of=process3,node distance=2cm] {Compute density $\mathbf{n}^{m+1}$, current density $\mathbf{J}^{m+1}$ and \\ momentum $\mathbf{S}^{m+1}$ by \eqref{eq:moment0-approx}--\eqref{eq:moment2-approx}};
		
		\draw [arrow] (input1) -- (process1);
		\draw [arrow] (process1) -- (process2);
		\draw [arrow] (process2) -- (process3);
		\draw [arrow] (process3) -- (stop);
	\end{tikzpicture}
	\caption{Flowchart of our AP-PIC method for solving the VPME system \eqref{RF_VPME1}--\eqref{RF_VPME2}.}
	\label{fig:flowchart}
\end{figure}
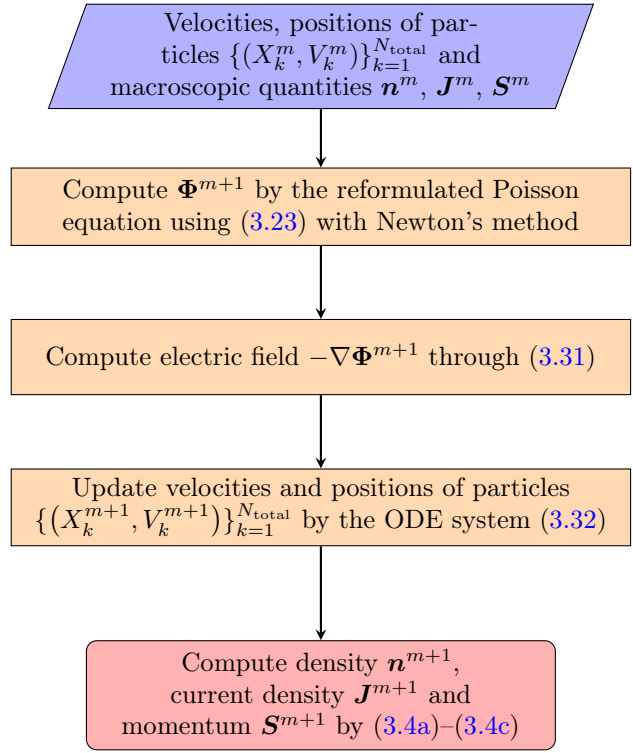

\begin{rmrk}
	\label{remark: penalty}
	Here, we provide an alternative approach to handle the nonlinear term $e^{\phi^{m+1}}$ in \eqref{RF-Poisson-time-discretation-2}, which does not require the Newton's iteration method. The idea of constructing a penalty term is inspired from \cite{lehman2022vlasov}. Consider  
	\begin{equation}
		\label{eq:penalty}
		e^{\phi} = \left[e^{\phi} - G(\phi)\right] + G(\phi). 
	\end{equation}
	To determine $G(\phi)$, a simple Taylor series expansion is performed around $\phi=0$ and we let 
	$G(\phi) = 1 + \phi$, then the nonlinear term $e^{\phi^{m+1}}$ becomes 
	\begin{equation}
		\label{eq:exp_phi_penalty}
		e^{\phi^{m+1}} \approx \left[e^{\phi^{m}} - (1+\phi^{m})\right] + (1+\phi^{m+1}).
	\end{equation}
	The discretized scheme for the reformulated Poisson equation \eqref{Full-scheme} turns out to be 
	\begin{equation}
		\label{eq:New-Full-scheme}
		\begin{array}{l}
			-\e^2 \left(\Delta\phi^{m+1}\right)_{j}- (\Delta t)^2 \Big(D^{-}\left(n_{i}^{m}D^{+} \phi^{m+1}\right)\Big)_{j} + \phi^{m+1}_{j}  \\[6pt]
			= (\Delta t)^2 \left(\Delta S_{i}^{m}\right)_{j}
			+  (n_i)_{j}^m -D^{-}(n_{i}u_{i})_j^m + (\phi^{m}_{j} - e^{\phi^{m}_{j}}), 
		\end{array}
	\end{equation}
	which can be solved explicitly. We show the scheme in a vector form: 
	\begin{equation}
		\label{eq:RF_Poisson_explicit}
		(A^{m} + I_{N_{h}\times N_{h}})\mathbf{\Phi}^{m+1} = \mathbf{b}^{m} + (\mathbf{\Phi}^{m} - e^{\mathbf{\Phi^{m}}}).
	\end{equation} 
	where $I$ is the unit matrix and $\mathbf{b}^m$ is defined in \eqref{Full-scheme-right}. As $\e\to 0$, the scheme \eqref{eq:RF_Poisson_explicit} automatically becomes a consistent discretization of the quasineutral limit equation \eqref{Quasi-Limit2}, if one applying the same approximation for $e^{\phi^{m+1}}$ in the limit system. Therefore, this new explicit scheme also satisfies the asymptotic-preserving property. In the first numerical test, we compute the reformulated Poisson by this explicit scheme \eqref{eq:RF_Poisson_explicit} and compare it with the implicit method using the Newton's iteration. This explicit scheme has shown advantages since the iteration is avoided, yet satisfies the AP property. To further explore this explicit scheme for the reformulated Poisson in the PIC framework and build higher-order AP method remains a future task. 
\end{rmrk}
\section{A bi-fidelity method for model with random parameters}
\label{sec:UQ}
In this section, we briefly review the bi-fidelity method. Assume the expensive high-fidelity model and the cheap low-fidelity model are available to generate the high-fidelity solution $\mathbf{U}^H(z)$ and the low-fidelity solution $\mathbf{U}^L(z)$ respectively, for any given parameter $z$. We refer to \cite{Narayan2014, ZNX2014} for more details and summarize the bi-fidelity method as follows: 
\begin{lgrthm}
	\label{AG}
\textbf{A bi-fidelity approximation}
	\begin{algorithmic}[1]
		\REQUIRE
		\STATE{Select a sample set $\Gamma = \{z_1, z_2, \hdots, z_{M}\}\subset I_z$, and run the low-fidelity model $\mathbf{U}^L(z_j)$ for each $z_j \in \Gamma$. }
		\STATE{Select $r$ ``important'' points from $\Gamma$ and denote it by $\gamma=\{z_1, \cdots z_{r} \} \subset\Gamma$. Construct the low-fidelity approximation space $\mathscr{U}^L(\gamma)$.}
		\STATE{Run high-fidelity simulation at each point in the selected set $\gamma$. Construct the high-fidelity approximation space $\mathscr{U}^H(\gamma)$.}
		\ENSURE
		\STATE{For any $z$, compute the low-fidelity solution $U^L(z)$ and the corresponding low-fidelity coefficients: }
		\begin{equation}\label{C-N}
			\mathbf{U}^L(z) \approx \mathcal{P}_{\mathscr{U}^L\left(\gamma\right)}\left[\mathbf{U}^L(z)\right]=\sum_{k=1}^r c_k^L(z) \mathbf{U}^L\left(z_k\right). \end{equation}
		\STATE{Construct the bi-fidelity approximation by applying the same approximation rule as in the low-fidelity model:}
		\begin{equation}\label{UB} \mathbf{U}^B(z) = \sum_{k=1}^r c_k^L(z) \mathbf{U}^H(z_k). \end{equation}
	\end{algorithmic}
\end{lgrthm}
Our sought bi-fidelity approximation of $u^H$ is constructed by equation \eqref{UB}. 
In the offline stage, one employs the cheap low-fidelity model to select the most important parameter points $\gamma=\{z_1, \cdots, z_r\}$ by the greedy procedure \cite{DeVore}. 
In the online stage, for any sample point $z\in I_z$, we project the low-fidelity solution $\mathbf{U}^L(z)$ onto the low-fidelity approximation space $\mathscr{U}^L(\gamma)$ by using \eqref{C-N}, where $\mathcal P_{\mathscr{U}^L(\gamma)}$ is the projection operator onto $\mathscr{U}^L(\gamma)$ with  coefficients $c_k$ computed by the Galerkin approach: 
$$ {\bf G}^L {\bf c} = {\bf f}, \, \text{~with~} 
{\bf f} = (f_k)_{1\leq k\leq r} \,  \text{~and~} 
f_k = \langle \mathbf{U}^L(z), \mathbf{U}^L(z_{k})\rangle^L, $$
where ${\bf G}^L$ is the Gramian matrix for $\mathbf{U}^L(\gamma)$. 

\vspace{2mm}

\noindent\textbf{High-fidelty solver:} We desire to solve the VPME system with multi-dimensional uncertain parameters, thus in our bi-fidelity method, the VPME system \eqref{VP-semi} is considered as the high-fidelity model. We employ the AP-PIC method studied in Subsection \ref{subsec:scheme} as our high-fidelity solver. 

\noindent\textbf{Low-fidelity solver:} We choose the macroscopic Euler system \eqref{Euler-semi} that characterizes a simpler physics as our low-fidelity model in the bi-fidelity method. We use the finite difference method (FDM) to discretize the spatial space and forward Euler scheme for the time discretization. Consider the conservation form of \eqref{Euler-semi} for one-dimensional space and velocity: 
\begin{equation}
	\label{Euler-semi-conserved}
	\frac{\partial\mathbf{U}}{\partial t} +\frac{\partial\mathbf{F}(\mathbf{U})}{\partial x}=0,\quad \mathbf{U} = \begin{bmatrix}
		n   \\nu
	\end{bmatrix},\quad \mathbf{F}(\mathbf{U})=\begin{bmatrix}
		nu   \\ nu^{2} + n
	\end{bmatrix}.
\end{equation}
where $n = n(x,t),\: u = u(x,t)$ are the ion density and velocity given in \eqref{ion-density-velocity} in the computational domain $\left[x_{L},x_{R}\right]$. 
Let $h$ be a uniform space step and $\tau$ be a uniform time step, then $x_{j} = jh,\:j=0,1,2,\dots,N_{l}$ with $x_{L} = x_{0},\:x_{R}=x_{N_{l}}$ and $t^{n} = n\tau,\: n=1,2,\dots,M$. 

In our flux construction, Lax-Friedrichs scheme--as one of the most popular finite difference methods for solving hyperbolic partial differential equations \cite{leveque1992numerical, thomas2013numerical, toro2013riemann} is used. 
Let $\mathbf{U}^{n}_{j} = \mathbf{U}(x_{j}, t^{n})$ be the approximated values of $\mathbf{U}$ at position $x_{j}$ and time $t^{n}$. The time discretization for \eqref{Euler-semi-conserved} becomes
\begin{equation}
	\label{Lax-scheme}
	\mathbf{U}^{n+1}_{j} = \mathbf{U}^{n}_{j} - \frac{\tau}{h}(\mathbf{F}^{n}_{j+\frac{1}{2}} -\mathbf{F}^{n}_{j-\frac{1}{2}}), 
\end{equation}
where $\mathbf{F}^{n}_{j+\frac{1}{2}}$ is given by 
\begin{equation}
	\label{Lax-flux}
	\mathbf{F}^{n}_{j+\frac{1}{2}} = \frac{1}{2}\left(\mathbf{F}(\mathbf{U}^{n}_{j+1}\right) + \mathbf{F}(\mathbf{U}^{n}_{j})) - \frac{h}{2\tau}(\mathbf{U}^{n}_{j+1}-\mathbf{U}^{n}_{j}). 
\end{equation}
Then the full scheme turns to be 
\begin{equation}
	\label{eq:Full-Lax-scheme}
	\mathbf{U}^{n+1}_{j} = \frac{\mathbf{U}^{n}_{j+1} + \mathbf{U}^{n}_{j-1}}{2}- \frac{\tau}{2h}(\mathbf{F}(\mathbf{U}^{n}_{j+1}) -\mathbf{F}(\mathbf{U}^{n}_{j-1})). 
\end{equation}
For the Euler system, we consider periodic boundary conditions 
\begin{equation}
	\label{eq: bc_for_Euler}
	n(x_{L},t) = n(x_{R}, t),\qquad u(x_{L},t) = u(x_{R},t). 
\end{equation}
\section{Numerical examples}
\label{sec:Num}
In our numerical experiments, we first study the deterministic VPME system, then extend to the problems with random uncertainties and validate the efficiency of our bi-fidelity method. For all the following experiments, the regularization parameter $\eta$ in \eqref{eq:delta-v-approx} and the error tolerance $\eta_{1}$ in the Newton's method \eqref{eq:Newton-iter-condition} for solving the reformulated Poisson are set as
$ \eta = 1, \, \eta_{1} = 1\times 10^{-8}$. 
We denote the total number of particles for a single PIC run as $N_{\text{tol}}$, and the notation $N_c$ represents the number of independent numerical simulations that are performed under the fixed setting and numerical parameters. The purpose of running the same test $N_c$ times 
is to reduce the statistical noise inherent in the PIC method through the ensemble averaging technique.
For notational convenience, we name the reformulated VPME system \eqref{RF_VPME1}--\eqref{RF_VPME2} as the RF-VPME system, and its quasineutral limit \eqref{Quasi-Limit1}--\eqref{Quasi-Limit2} as the Quasi-VP system.  

\subsection{The RF-VPME system: an explicit scheme}
\label{sec:penalty}
\begin{figure}[!hptb]
	\centering
	\subfigure[]{
		\includegraphics[width = 0.4\textwidth]{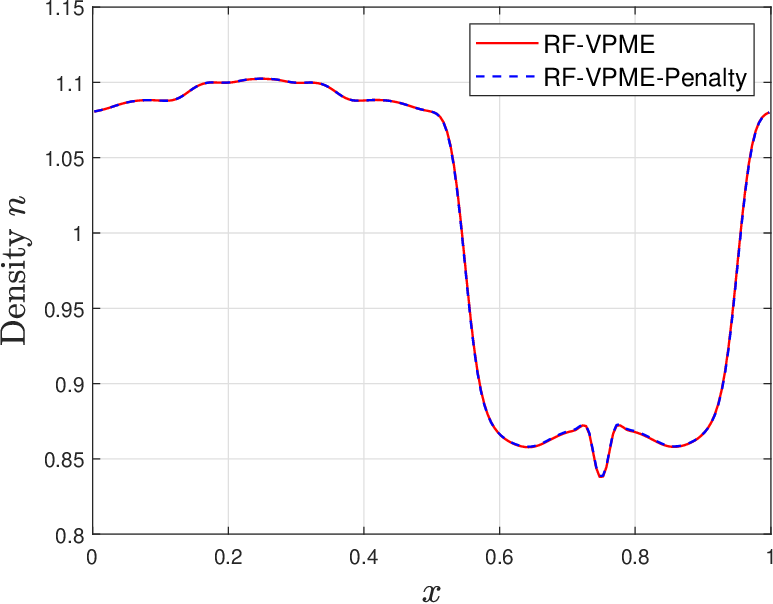}
	}
	\subfigure[]{               
		\includegraphics[width = 0.4\textwidth]{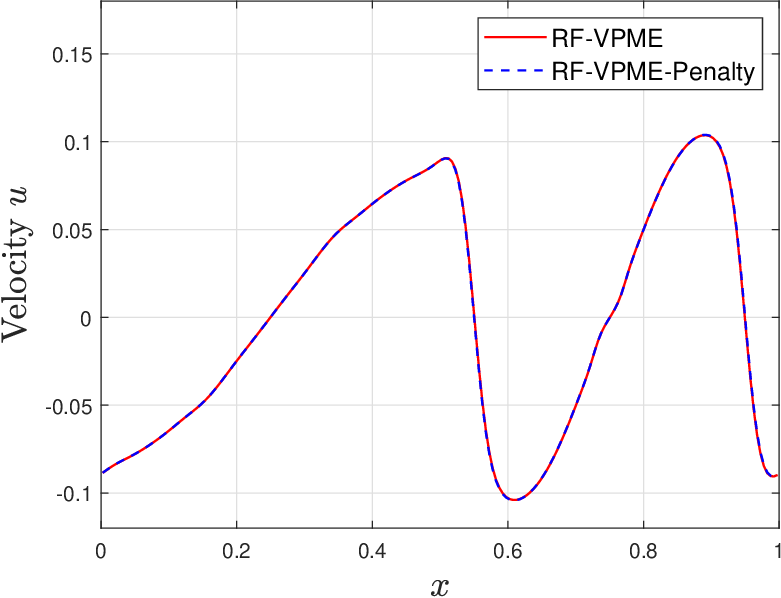}
	}       
	\caption{(Two schemes for solving the reformulated Poisson, $\e = 0.01$). Density $n$ (left) and velocity $u$ (right) at $t=2$ by an implicit scheme \eqref{Full-scheme} with Newton's iteration (red) and an explicit scheme based on penalty \eqref{eq:RF_Poisson_explicit} (blue).}    
	\label{fig:poisson_compare_ep=0.01}
\end{figure}

Corresponding to Remark \ref{remark: penalty} based on the penalty, in this first test we compute the reformulated Poisson by an explicit scheme \eqref{eq:RF_Poisson_explicit} and compare the solution with that obtained from the Newton's method. For simplicity, a finite difference scheme \cite{leveque1992numerical, toro2013riemann} in the Eulerian framework with discrete velocity method (DVM) \cite{cabannes1980discrete, sharipov2015rarefied} for solving the Vlasov equation \eqref{eq:VPME-Vlasov} is adopted.
Let the initial condition $f_0$ be 
\begin{equation}
	\label{eq:ini_penalty_poisson} 
	f(0, x, v) = n_0(x)\delta(v-u_0(x)),
\end{equation}
with 
\begin{equation}
	\label{eq:poisson-test}
	n_{0}(x) = 1 + 0.5\sin{(2\pi x)},\qquad u_{0}(x) = 0.4\cos{(2\pi x)}, \qquad x \in D=[0, 1].
\end{equation}
Set the final time as $t=2$ and $\e=0.01$ in the model. Let $v\in [-4,4]$ with uniform mesh and the number of grid points be $N_v=800$, the spatial mesh size be $N_x=200$ and time step be $\Delta t=0.2\Delta x$. 
The density $n$ and velocity $u$ are shown in Figure \ref{fig:poisson_compare_ep=0.01}, where the blue line represents solution of the explicit scheme \eqref{eq:RF_Poisson_explicit} based on building a penalty for the nonlinear term $e^{\phi^{m+1}}$, and the red line stands for the implicit scheme using Newton's iteration \eqref{Newton-method-update}. The solutions of both methods match well, and both schemes are stable for long time. 
In our future work, we will incorporate this explicit scheme with the PIC method for solving the Vlasov equation \eqref{eq:VPME-Vlasov} and develop higher-order AP schemes. 

\subsection{The RF-VPME system: deterministic problem}
\label{sec:AP}
In this section, the AP property for the RF-VPME system is tested from three different perspectives. 
\begin{equation}
	\label{eq:ini_AP} 
	f(0, x, v) = n_0(x)\delta(v-u_0(x)),
\end{equation}
with $n_0(x)$ and $u_0(x)$ defined by  
\begin{equation}
	\label{eq:AP-test}
	n_{0}(x) = 1 + 0.5\sin{(2\pi x)},\qquad u_{0}(x) = 0.4\cos{(2\pi x)}, \qquad x \in D=[0, 1].
\end{equation}

In this test, we examine the convergence in space  of our numerical scheme for the RF-VPME system. The total number of particles in the AP-PIC method is set as $N_{\rm total} = 1.6\times 10^{6}$, with $N_c = 300$ times of repetition in the PIC method to reduce the noise. The CFL number is chosen as ${\rm CFL} = 0.4$. We obtain the reference solution by using spatial mesh size $N_x=800$. 

\begin{figure}[!hptb]
	\centering
	\subfigure[density $n$ ($\varepsilon = 1$)]{
		\includegraphics[width = 0.30\textwidth]{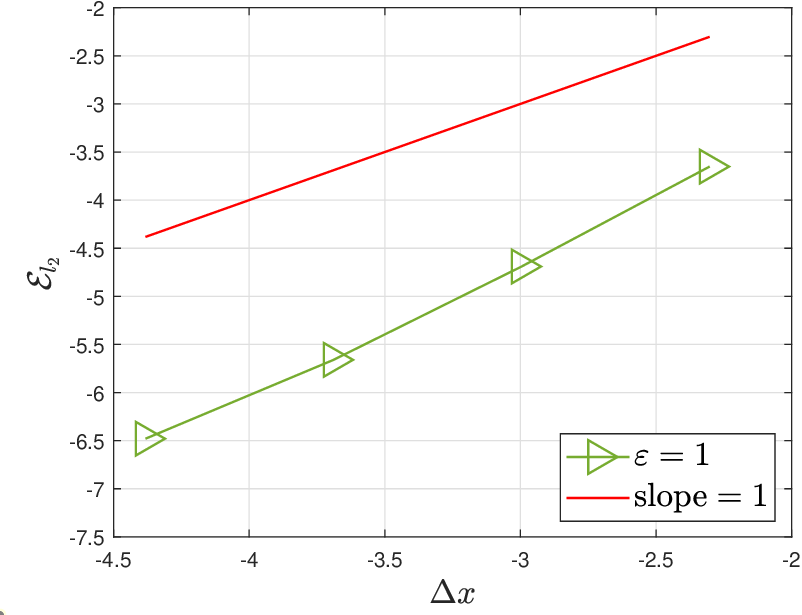}
	}
	\subfigure[velocity $u$ ($\varepsilon = 1$)]{               
		\includegraphics[width = 0.30\textwidth]{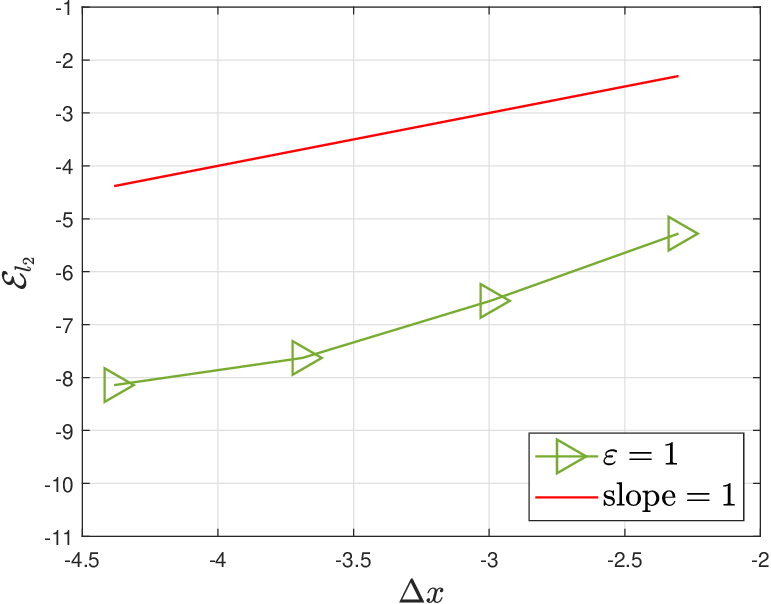}
	}       
	\subfigure[potential $\phi$ ($\varepsilon = 1)$]{
		\includegraphics[width = 0.30\textwidth]{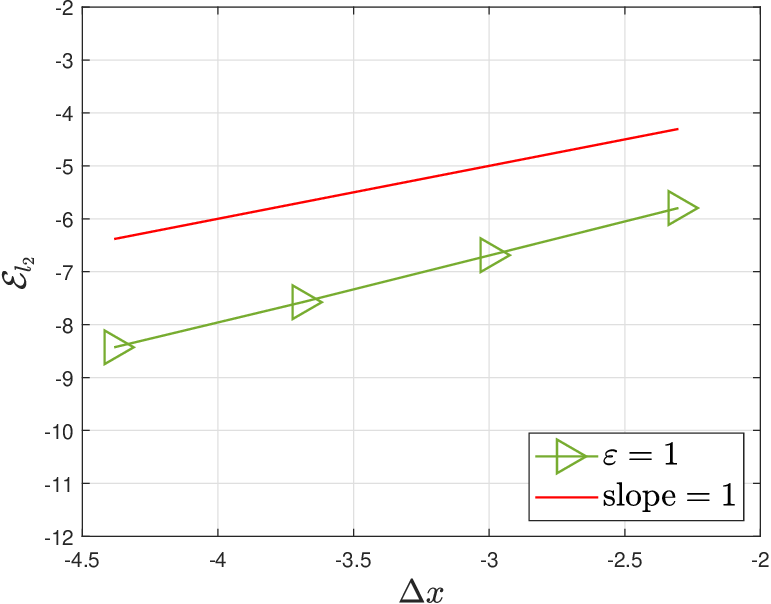}
	}
	\caption{(AP property for the RF-VPME in Subsection \ref{sec:AP}). The $l_{2}$ errors $\mathcal{E}_{l_2}$ of density $n$, velocity $u$ and potential $\phi$ for $\varepsilon = 1$ at $t = 0.2$. Here, the $x$-axis is $\log(\Delta x)$, and the $y$-axis is $\log(\mathcal{E}_{l_2})$. The mesh sizes in space are $N_x = 10, 20, 40$ and $50$. }    
	\label{fig:AP_conv_ep=1}
\end{figure}

\begin{figure}[!hptb]
	\centering
	\subfigure[density $n$]{
		\includegraphics[width = 0.30\textwidth]{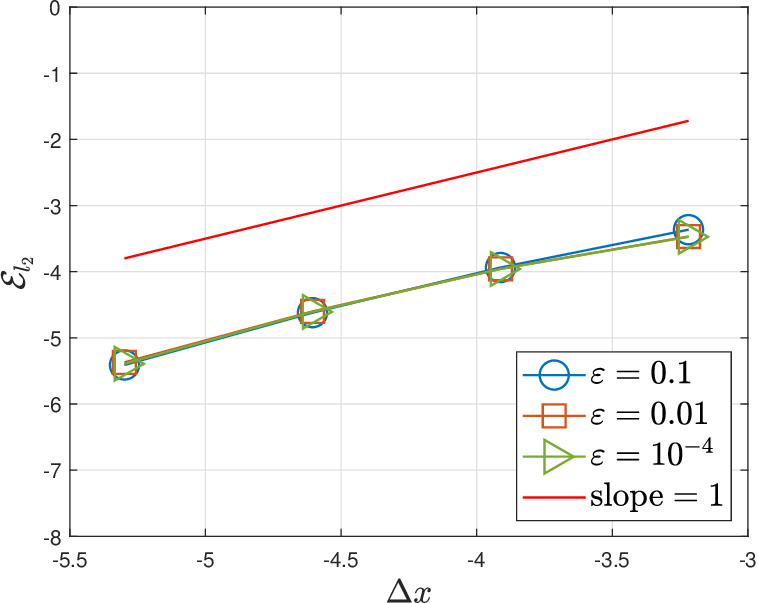}
	}
	\subfigure[velocity $u$]{               
		\includegraphics[width = 0.30\textwidth]{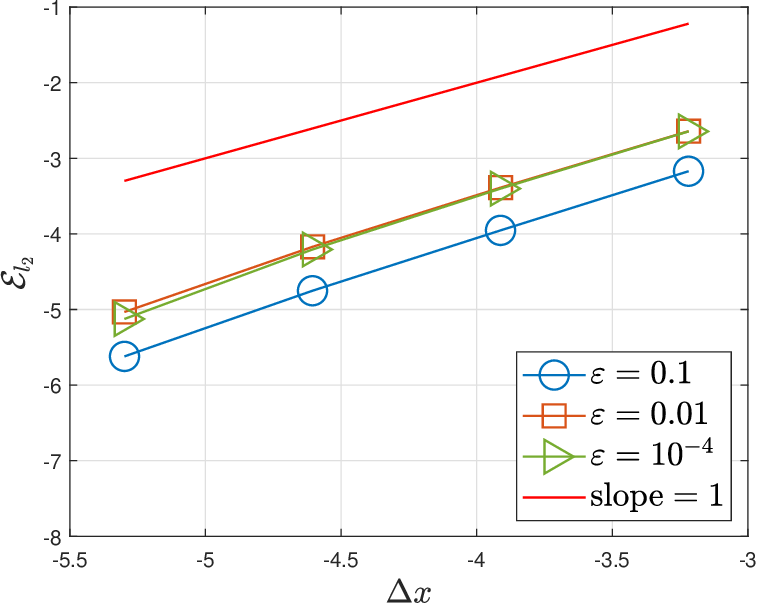}
	}       
	\subfigure[potential $\phi$ ]{
		\includegraphics[width = 0.30\textwidth]{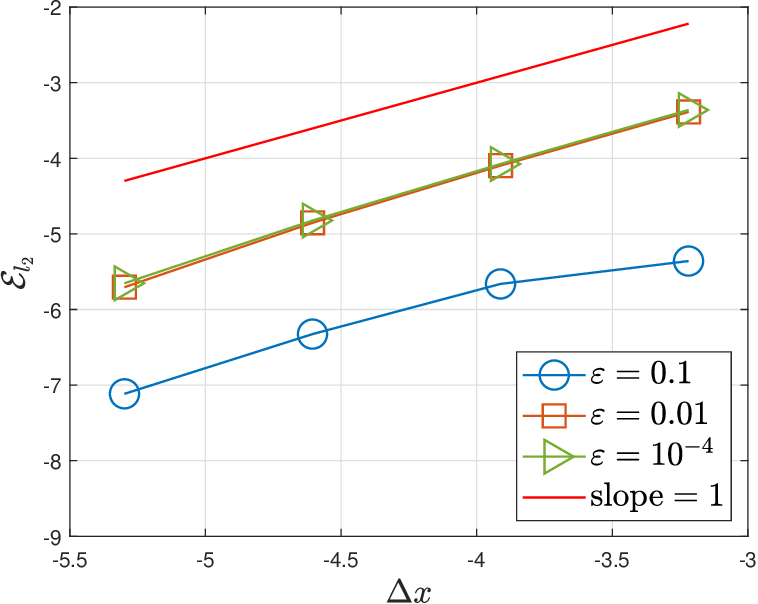}
	}
	\caption{(AP property for the RF-VPME in Subsection \ref{sec:AP}). The $l_{2}$ errors $\mathcal{E}_{l_2}$ of density $n$, velocity $u$ and potential $\phi$ for $\varepsilon = 0.1, 0.01$ and $10^{-4}$ at $t = 0.2$. Here, the $x$-axis is $\log(\Delta x)$, and the $y$-axis is $\log(\mathcal{E}_{l_2})$. The mesh sizes in space are $N_x = 25, 50, 100$ and $200$.} 
	\label{fig:AP_conv_ep_le_1}
\end{figure}

We particularly study the cases when $\varepsilon=1, 0.1, 0.01$ and $10^{-4}$. The final time is $t=0.2$. When $\varepsilon = 1$, the numerical solution is obtained by letting $N_x = 10, 20, 40$ and $50$. The $l_2$ errors $\mathcal{E}_{l_2}$ between numerical solutions and the reference solution for $\varepsilon = 1$ are shown in Figure \ref{fig:AP_conv_ep=1}, where the first-order convergence is clearly observed. 
For other $\varepsilon$, the $l_2$ errors $\mathcal{E}_{l_2}$ between the numerical solutions using $N_x = 25, 50, 100$ and $200$ and the reference solution is illustrated in Figure \ref{fig:AP_conv_ep_le_1}. For all these $\varepsilon$, the first-order convergence is obviously detected, indicating the efficiency of our AP-PIC scheme and a {\it uniform} first-order accuracy in space for different orders of $\varepsilon$. Therefore, our numerical scheme satisfies the AP property.

In the second test, we study the long-time behavior of the solutions to the RF-VPME system by comparing it with the solutions obtained by solving the quasineutral limit \eqref{Quasi-Limit1}-\eqref{Quasi-Limit2}, namely the Quasi-VP system. We consider the cases when $\e = 0.1,0.01, 10^{-4}$ and plot the solutions from the initial condition to $t=15$. 
At fixed time, the relative error between the solutions computed by the RF-VPME and the Quasi-VP system is defined as  
\begin{equation}
	\label{relative-error}
	\mathcal{E}_{r}(g) = \frac{\sqrt{\sum_{i=1}^{N_{x}}(g_{i}-q_{i})^2}}{\sqrt{\sum_{i=1}^{N_{x}}q_{i}^{2}}}.
\end{equation}
where $g$ represents the quantities (density $n$, velocity $u$ or potential $\phi$) solved by the RF-VPME system, and $q$ represents the corresponding solutions obtained by the Quasi-VP system. 

\begin{figure}[!hptb]
	\centering
	\subfigure[density $n$]{
		\includegraphics[width = 0.30\textwidth]{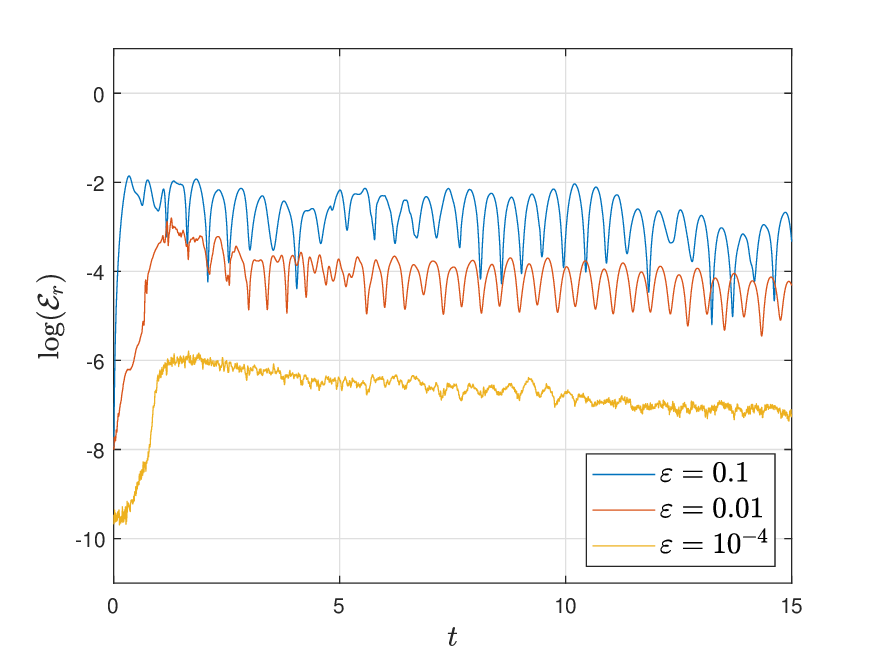}
	}
	\subfigure[velocity $u$]{               
		\includegraphics[width = 0.30\textwidth]{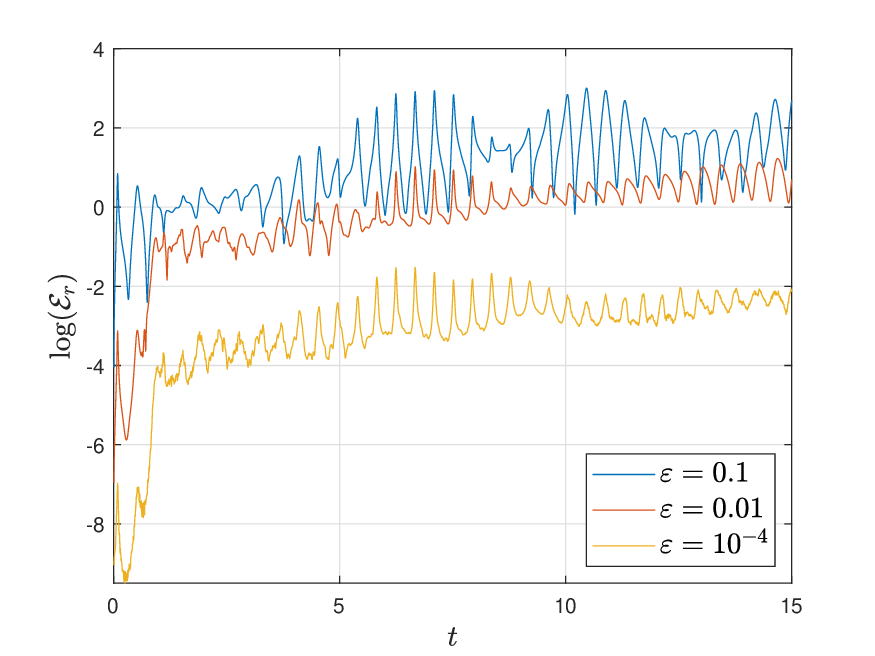}
	}       
	\subfigure[potential $\phi$ ]{
		\includegraphics[width = 0.30\textwidth]{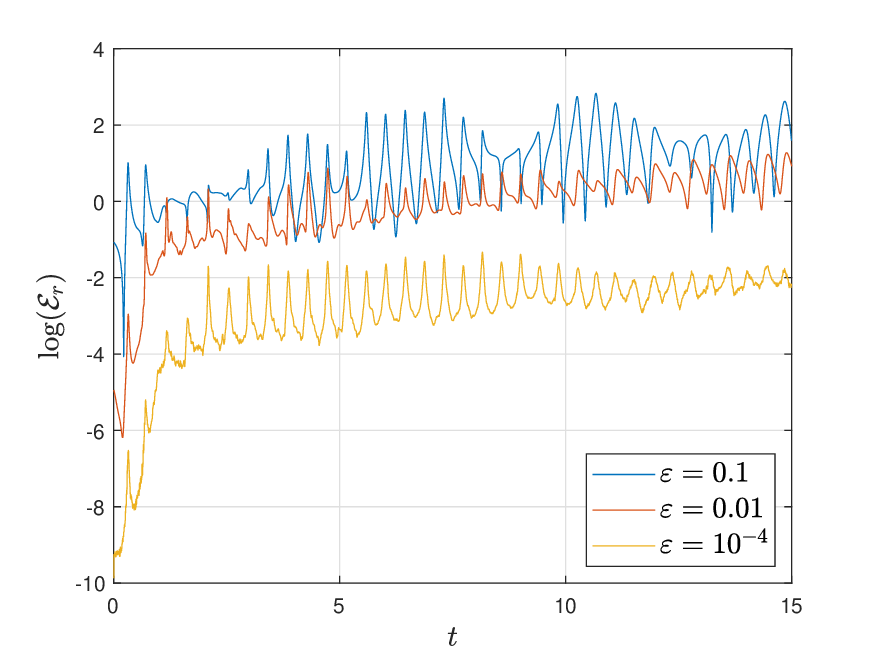}
	}
	\caption{(Long time evolution of RF-VPME in Subsection \ref{sec:AP}). The relative errors $\mathcal{E}_{r}$ of density $n$, velocity $u$ and potential $\phi$ with $\varepsilon = 0.1, 0.01$ and $10^{-4}$ between RF-VPME and Quasi-VP models, from $t = 0$ to $t = 15$. Here, the $x$-axis is time $t$, and the $y$-axis is $\log(\mathcal{E}_{r})$. The mesh size in space is $N_x=100$. } 
	\label{fig:relative-error}
\end{figure}

In the AP-PIC method, the number of total particles is set as $N_{\rm total} = 1.6\times 10^{6}$, with $N_c = 200$ times of repetition for reduction of noise. Let the CFL number be ${\rm CFL} = 0.2$. 
We use $N_x=100$ to numerically solve the Quasi-VP and RF-VPME systems. 
The relative errors $\mathcal{E}_{r}$ between the solutions of two systems for $\varepsilon = 0.1,0.01,10^{-4}$ are shown in Figure \ref{fig:relative-error}. One can observe that the smaller $\varepsilon$ is, the smaller the relative errors between the two systems are, which satisfies our expectation since the difference between the RF-VPME and Quasi-VP should be decreasing as $\e$ approaches zero. In Table \ref{table:relative_error_t=5}, we show the relative errors of $n$, $u$ and $\phi$ at $t=5$ for $\e=0.1,0.01,10^{-4}$, where the same trend of decreasing errors is detected. Therefore, for long-time simulations, our scheme for solving the RF-VPME system is stable and the solution converges to that of the quasineutral limiting equation as $\e$ becomes small. In addition, from Table \ref{table:relative_error_t=5} we notice that the space discretizations do not need to resolve the Debye length as it stays fixed ($N_x=100)$ for different ranges of $\e$, which validates again the AP property of our numerical scheme. 
\begin{table}[htbp]
	\centering
	\def\arraystretch{1.5}
		\begin{tabular}{c|ccc}
			\multicolumn{4}{c}{{\bf relative error $\mathcal{E}_{r}(g)$ at $t = 5$}}\\
			\hline
			$\varepsilon$ & $\mathcal{E}_{r}(n)$ & $\mathcal{E}_{r}(u)$ & $\mathcal{E}_{r}(\phi)$\\
			\hline
			$1\times 10^{-1}$    & 0.0547 & 1.2688 & 1.1648 \\
			\hline
			$1\times 10^{-2}$    & 0.0147 & 0.4886 & 0.2755 \\
			\hline
			$1\times 10^{-4}$  & 0.0025 & 0.0209 & 0.0407\\  
			\hline 
		\end{tabular}
        \vspace{4mm}
	\caption{(Relative errors between RF-VPME and Quasi-VP system at $t=5$ in Subsection \ref{sec:AP}). The relative errors of density $\mathcal{E}_{r}(n)$, velocity $\mathcal{E}_{r}(u)$ and potential $\mathcal{E}_{r}(\phi)$ for $\e = 0.1, 0.01, 10^{-4}$ at $t=5$. The mesh size in space is $N_x=100$.
	}
	\label{table:relative_error_t=5}
\end{table}

\begin{figure}[!hptb]
	\centering
	\includegraphics[width=0.5\linewidth]{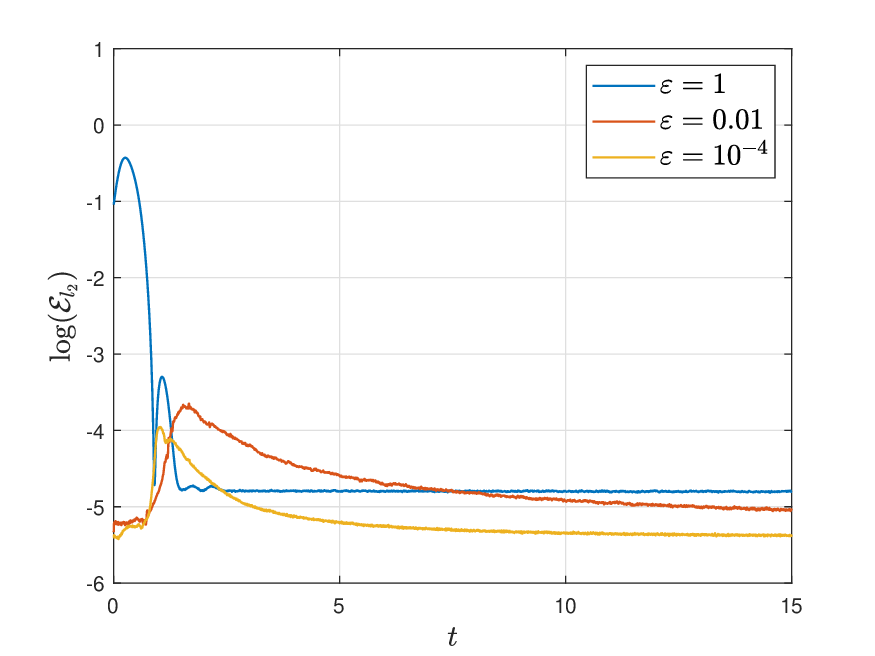}
	\caption{(Long time evolution for the errors between $n_{i}$ and $e^{\phi}$ in Subsection \ref{sec:AP}). The $l_{2}$ errors $\mathcal{E}_{l_2}$ of $(n_{i} - e^{\phi})$ for $\varepsilon = 1, 0.01, 10^{-4}$ from $t = 0$ to $t = 15$. Here, the $x$-axis is time $t$, and the $y$-axis is $\log(\mathcal{E})$. The mesh size in space is $N_x = 100$.}
	\label{fig:error-rho-expphi}
\end{figure}

In the last test of this section, we investigate the $l_2$ errors between $n_i$ and $e^{\phi}$ as time evolves and for different $\e$ in the model. In the AP-PIC method, we set $N_{\rm total} = 1.6\times 10^{6}$, $N_c = 200$, $N_x=100$, and the CFL number is ${\rm CFL} = 0.2$.
Three cases of $\e = 1, 0.01, 10^{-4}$ are considered.
In each simulation, we compute $l_2$ errors of \( (n_i - e^{\phi}) \) and obtain the error values by taking the average of these $N_c$ groups of errors. In Figure \ref{fig:error-rho-expphi}, we notice that as $\e$ becomes smaller, the errors between $n_i$ and $e^{\phi}$ decrease corresponding, especially at the longer time when errors saturate.  This is because that the equation \eqref{P2}, together with initial conditions \eqref{IC-phi} that we used in the AP-PIC scheme, is equivalent to the original Poisson equation \eqref{Poisson}, which leads to the relation $n_i = e^{\phi}$ when $\e$ vanishes. 

\subsection{The RF-VPME system: uncertain problem}

In this section, we extend to study the RF-VPME system with uncertain parameters and investigate the efficiency of our bi-fidelity (BF) approach as introduced in Subsection \ref{sec:UQ}. 
The high-fidelity (HF) model is the RF-VPME system \eqref{RF_VPME1}--\eqref{RF_VPME2} 
solved by the AP-PIC method shown in Figure \ref{fig:flowchart}. The low-fidelity (LF) model is chosen as the Euler system \eqref{Euler-semi-conserved} computed by using the Lax-Friedrichs scheme. To evaluate the accuracy of our method, we define errors between the BF and HF solutions at time $t$: 
\begin{equation}
	\label{eq:bi-high-l2-error}
	\left\|U^{H}(t)-U^{B}(t)\right\|_{L^{2}\left(D \times I_{z}\right)} \approx \frac{1}{K} \sum_{i=1}^{K}\left\|U^{H}\left(\hat{z}_{i}, t\right)-U^{B}\left(\hat{z}_{i}, t\right)\right\|_{L^{2}(D)},
\end{equation}
where $D$ is the spatial domain, and the set of points $\{\hat{z}_{i}\}_{i=1}^{K} \subset I_{z}$ is sampled independently of $\Gamma$ set in Algorithm \ref{AG}. In the following tests, we let the training set $\Gamma$ be $M=1000$ randomly selected points on $I_z$ and examine the errors of bi-fidelity approximation using \eqref{eq:bi-high-l2-error} with $K=200$.

\subsubsection{Test 1: uncertain initial conditions}
\label{sec:uq_test1}
In this test, we consider the case where the initial data contains random parameters. Assume that
\begin{equation}
	\label{eq:ini_uq}
	f(0,x,v,\mathbf{z}) = n_{0}(x,\mathbf{z})\delta(v-u_{0}(x)),
\end{equation}
with $n_{0}(x,\mathbf{z})$ and $u_{0}(x)$ given by
\begin{equation}
	\label{eq:ini_new_uq1}
	\begin{aligned}
		n_{0}(x,\mathbf{z})  = & 1 + 0.5\cos{(2\pi x)} + \frac{1}{\pi^2}\cos(2\pi x)z_1 +  
		+\frac{1}{(2\pi)^2}\sin(2\pi x)z_2 \\
		& +  \frac{1}{(3\pi)^2}\cos(2\pi x)z_3 
		+\frac{1}{(4\pi)^2}\sin(2\pi x)z_4 
		+ \frac{1}{(5\pi)^2}\cos(2\pi x)z_5, \\[4pt] 
		u_{0}(x) &= 0.2\sin{(2\pi x)},
		\quad 
		x\in\left[0,1\right], 
	\end{aligned}
\end{equation}
where the random variable $\mathbf{z}= (z_1, \cdots, z_5)$ is 5-dimensional, with $\{z_i\}_{i=1}^5$ following the uniform distribution on $I_z=[-1,1]$. We study the model for different $\e$ and let final computational time $t = 0.5$. In the HF solver, we adopt the spatial mesh size $N_{h} = 400$, total number of particles $N_{tol} = 1.6\times10^{6}$ and repeatedly simulate $N_c=200$ times to reduce noise in the PIC method; the time step satisfies the stability condition and is set as
$ \Delta t_{h} = 0.4\frac{\Delta x}{\max\{u_0\}}$.   
In the LF solver, we take the spatial mesh size $N_{l} = 1000$,  
employ the Lax-Friedrichs scheme and let 
$ \Delta t_{l} = 0.1 \Delta x$.

\begin{figure}[!hptb]
	\centering
	\subfigure{
		\includegraphics[width = 0.4\textwidth]{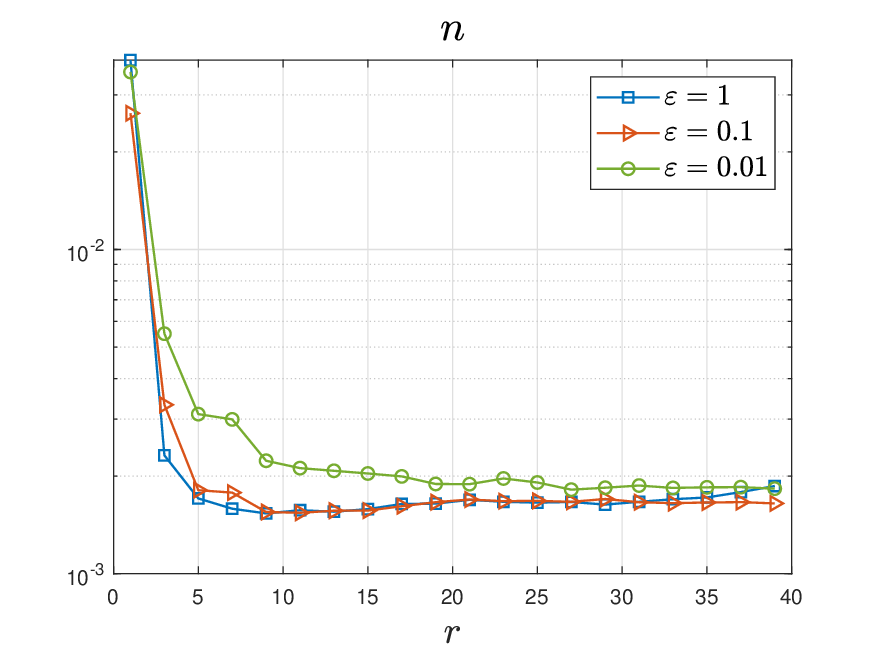}
	}
	\subfigure{               
		\includegraphics[width = 0.4\textwidth]{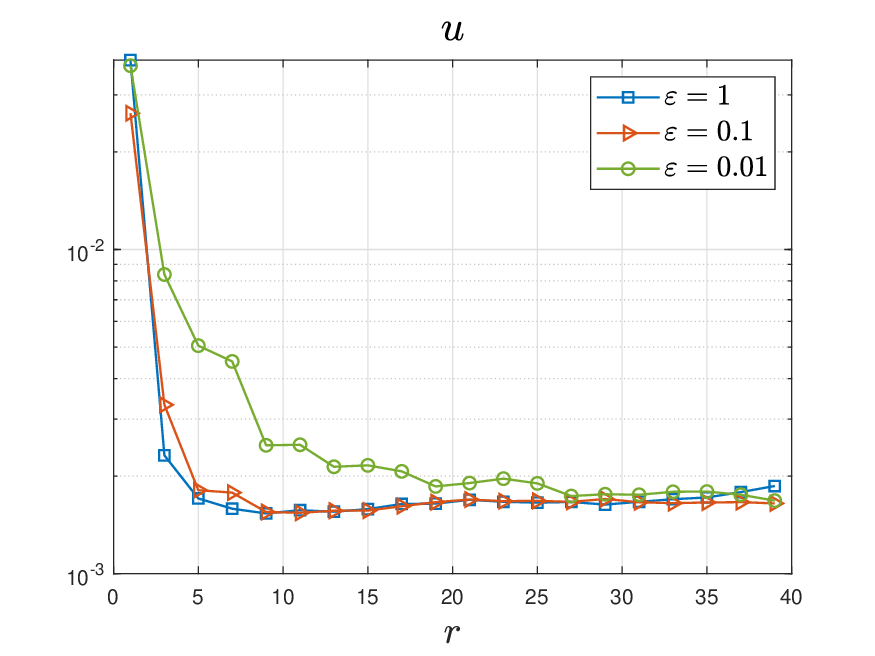}
	}       
	\caption{(UQ Test 1, different $\e$ in Subsection \ref{sec:uq_test1}) Errors of the bi-fidelity approximation for $n$ (left), $u$ (right) with respect to the number of high-fidelity simulation runs.}    
	\label{fig:uq_init2_error}
\end{figure}

\begin{figure}[!hptb]
	\centering
	\subfigure{
		\includegraphics[width = 0.4\textwidth]{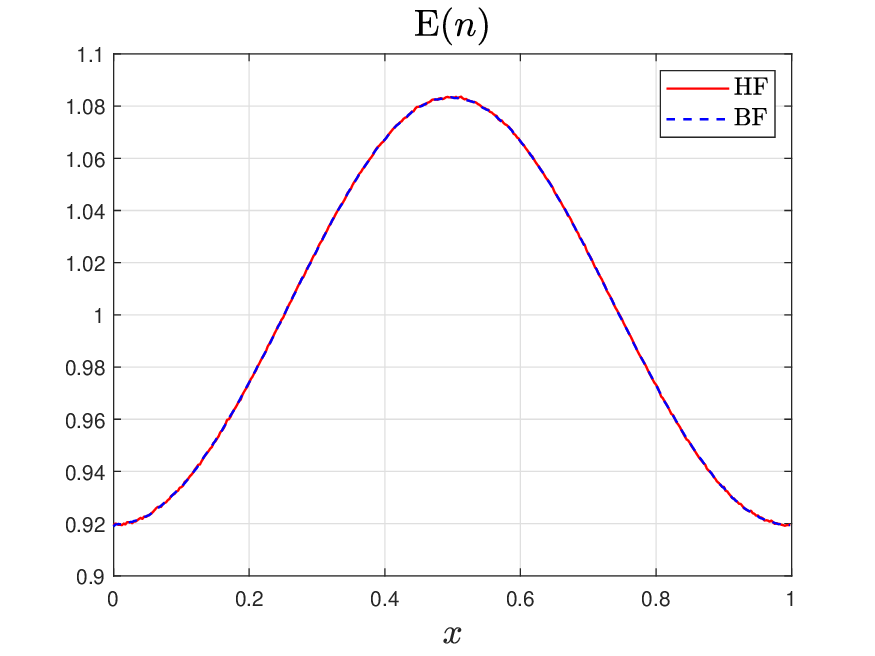}
	}
	\subfigure{               
		\includegraphics[width = 0.4\textwidth]{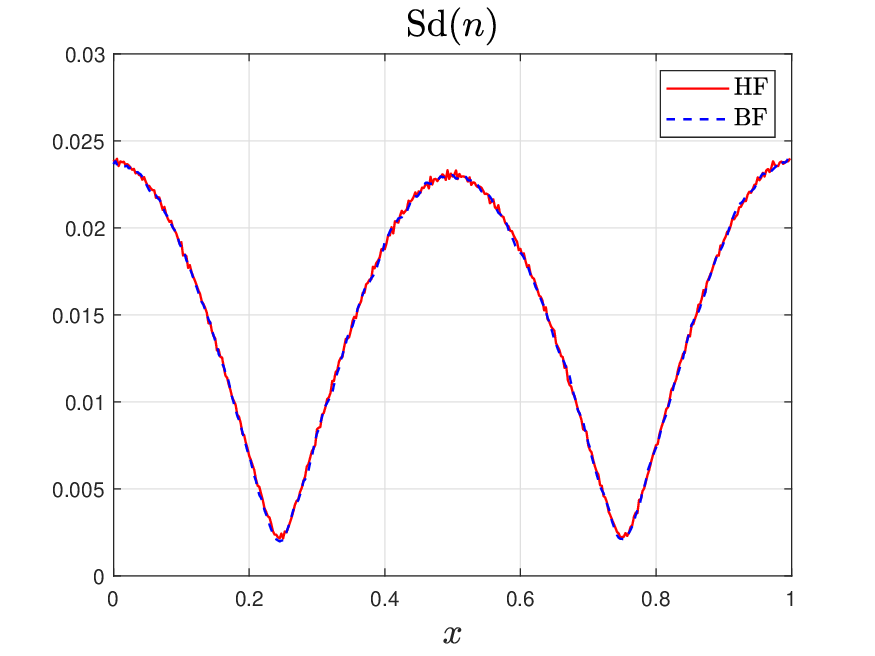}
	}       
	\subfigure{
		\includegraphics[width = 0.4\textwidth]{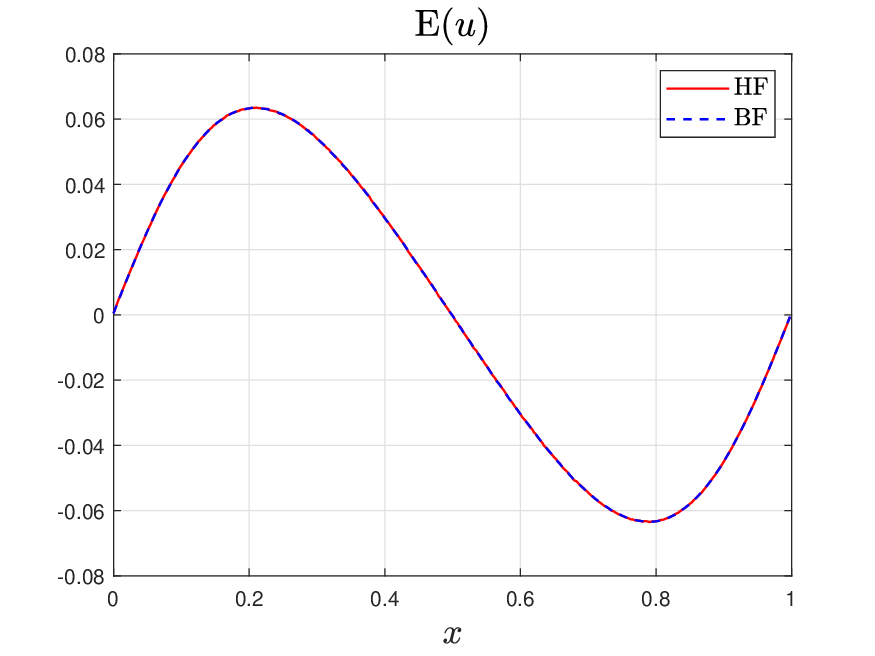}
	}
	\subfigure{               
		\includegraphics[width = 0.4\textwidth]{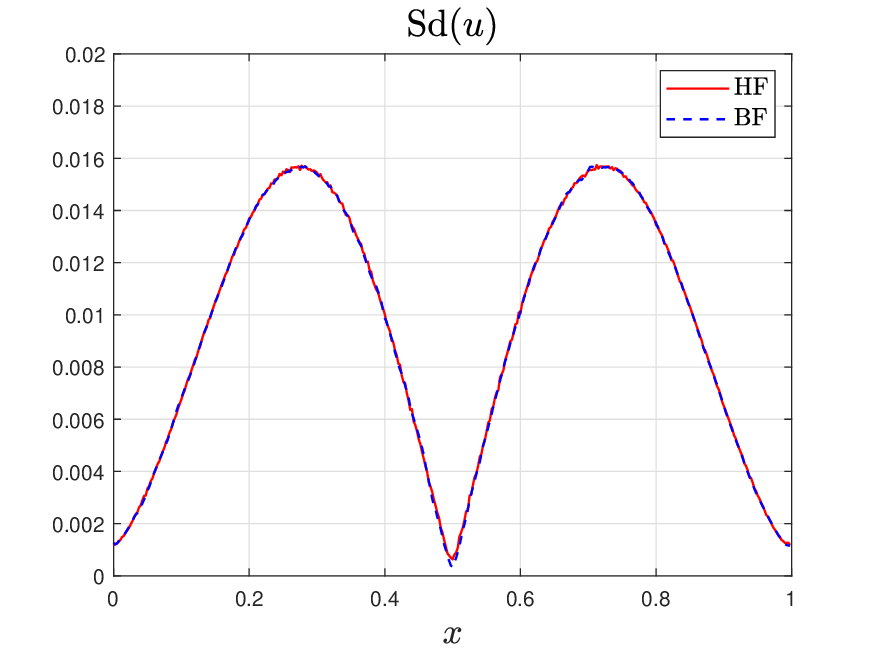}
	} 
	\caption{
		(UQ Test 1, $\e=1$ in Subsection\ref{sec:uq_test1}) Mean and standard deviation of $n$, $u$ of high-fidelity and bi-fidelity solutions by using $r=10$. }    
	\label{fig:uq_1e0_init2_mean_std}
\end{figure}

In Figure \ref{fig:uq_init2_error}, we consider the model under different regimes with $\e = 1, 0.1, 0.01$. This figure shows the errors of density $n$ and velocity $u$ between BF and HF solutions, as the number of HF simulation runs $r$ increases. A fast exponential decay of the errors is clearly observed in both plots, for different cases of $\e$. In addition, we can detect that the errors for these macroscopic quantities saturate as $r$ reaches about 10, at the accuracy level of $O(10^{-3})$. In Figure \ref{fig:uq_1e0_init2_mean_std}, we test that even for large $\e$ ($\e = 1$), the mean and standard deviation of BF approximations match really well with the HF solutions, by running the HF solver only 10 times. 

\begin{figure}[!hptb]
	\centering
	\subfigure[$\e = 1$, density $n$]{
		\includegraphics[width = 0.3\textwidth]{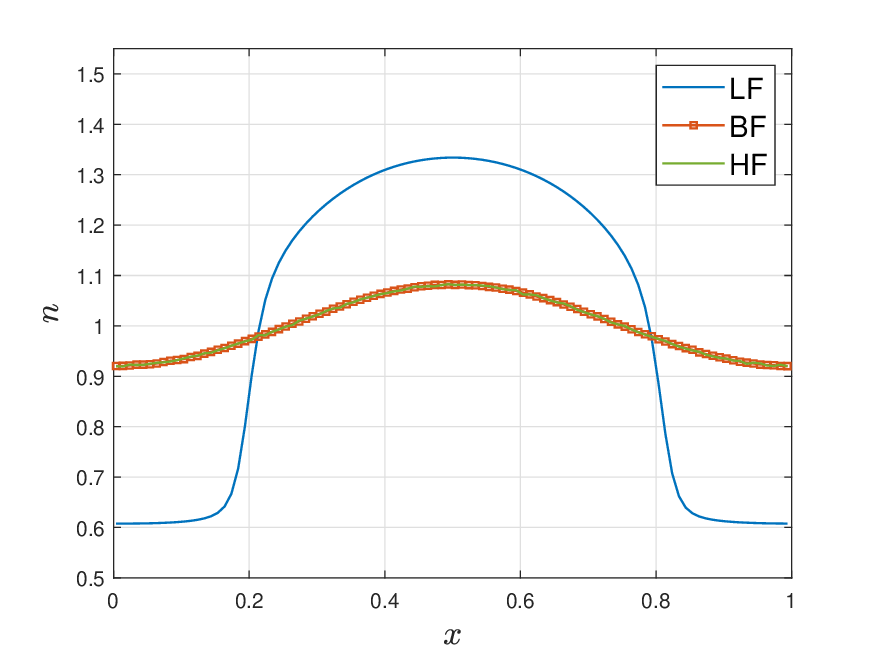}
	}
	\subfigure[$\e = 0.1$, density $n$]{
		\includegraphics[width = 0.3\textwidth]{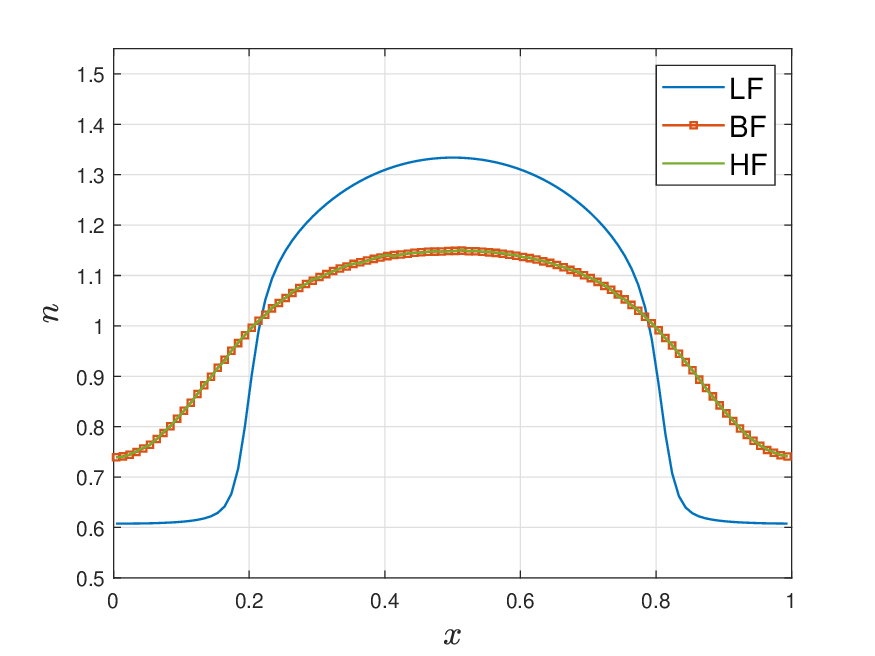}
	}
	\subfigure[$\e = 0.01$, density $n$]{
		\includegraphics[width = 0.3\textwidth]{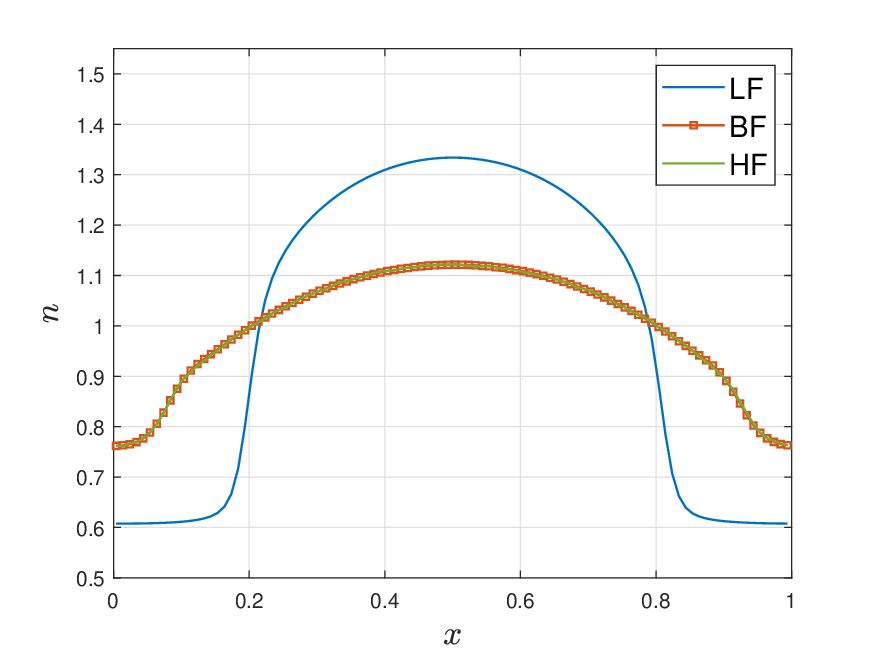}
	}
	\subfigure[$\e = 1$, velocity $u$]{               
		\includegraphics[width = 0.3\textwidth]{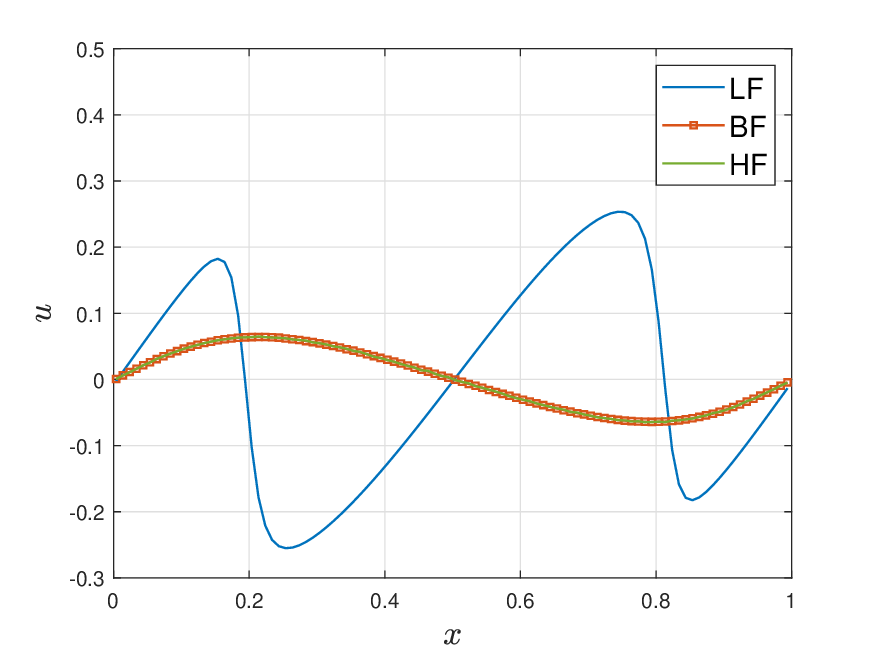}
	}     
	\subfigure[$\e = 0.1$, velocity $u$]{               
		\includegraphics[width = 0.3\textwidth]{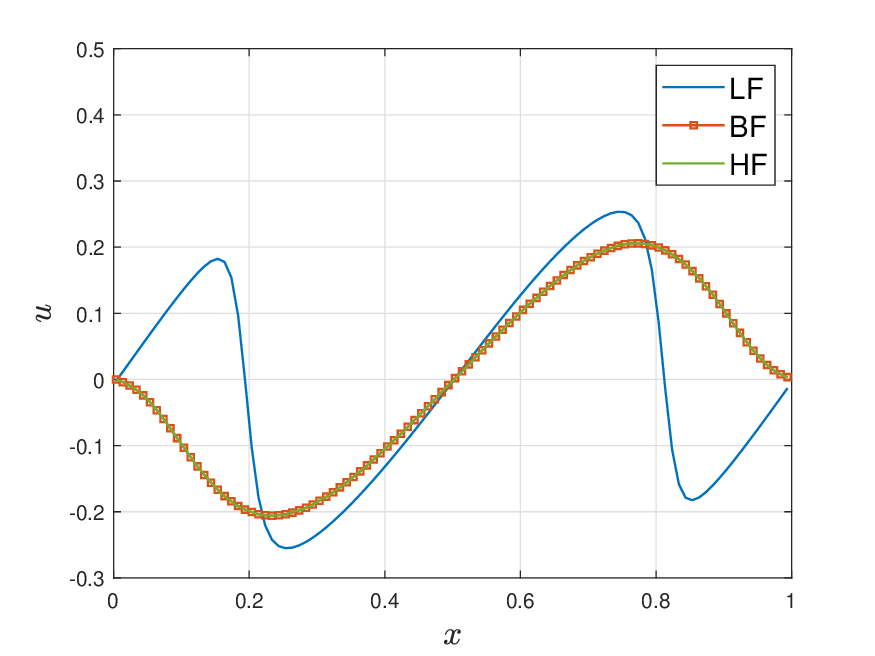}
	}
	\subfigure[$\e = 0.01$, velocity $u$]{               
		\includegraphics[width = 0.3\textwidth]{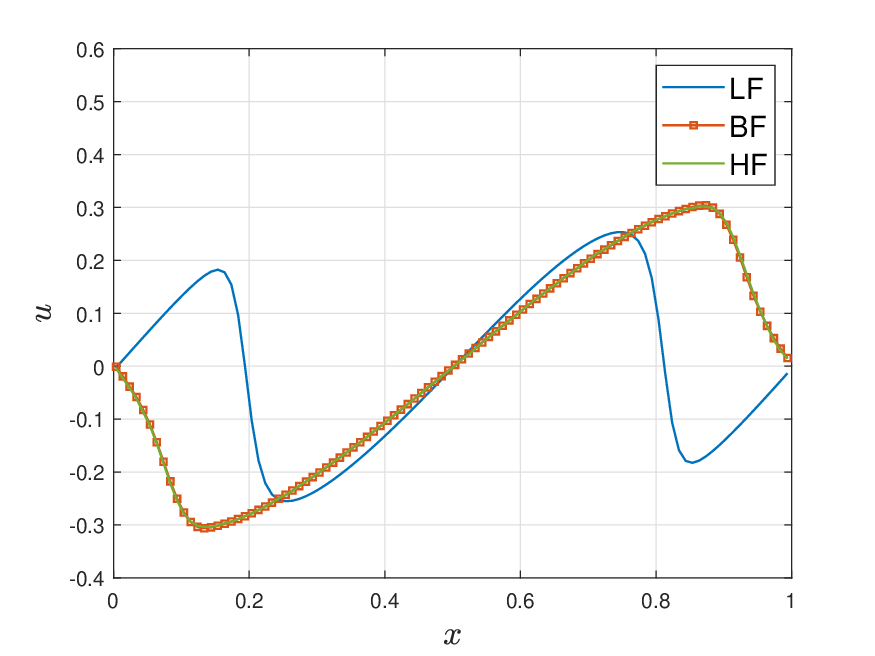}
	} 
	\caption{
		(UQ Test 1, different $\e$ in Subsection \ref{sec:uq_test1}) Comparison of low-, high- and bi-fidelity solutions at a fixed $z$ by using $r=10$ in the bi-fidelity approximations. }    
	\label{fig:uq_init2_fixed_z}
\end{figure}

In Figure \ref{fig:uq_init2_fixed_z}, we consider the cases $\e=1,0.1,0.01$ and compare the HF, LF and BF solutions of $n$ and $u$ at an arbitrarily selected $z$ point, by using $r=10$ in the BF approximation. It can be seen that the HF and BF solutions are consistent with each other, whereas the LF solutions display a relatively large difference at many spatial locations.  
Together with Figure \ref{fig:uq_1e0_init2_mean_std}, these results suggest that even though our LF Euler model may not be accurate in the physical space for each fixed $z$, it still can capture the behavior and characteristics of the HF model (VPME system) in the random space, thus demonstrating the effectiveness and accuracy of our proposed BF approach. 

\subsubsection{Test 2: uncertain initial data in the double-delta form}
\label{sec:uq_test2}
\begin{figure}[!hptb]
	\centering
	\subfigure{
		\includegraphics[width = 0.3\textwidth]{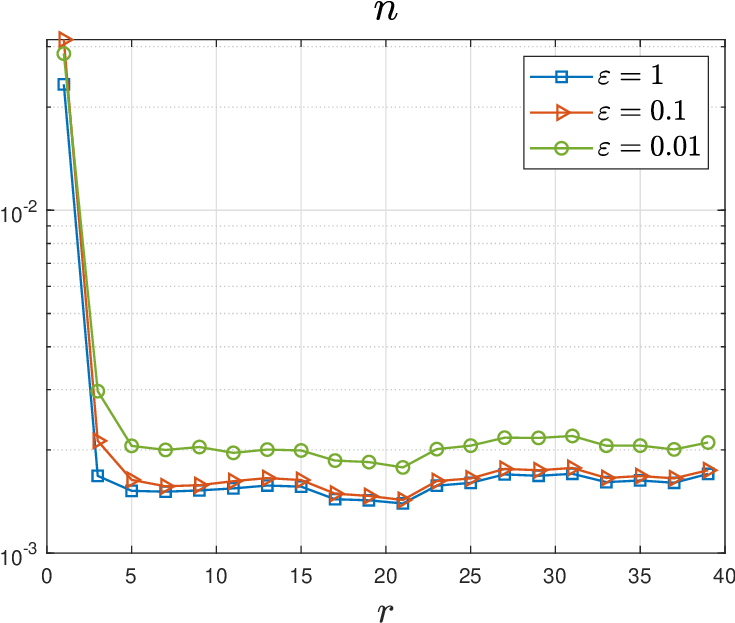}
	}
	\subfigure{               
		\includegraphics[width = 0.3\textwidth]{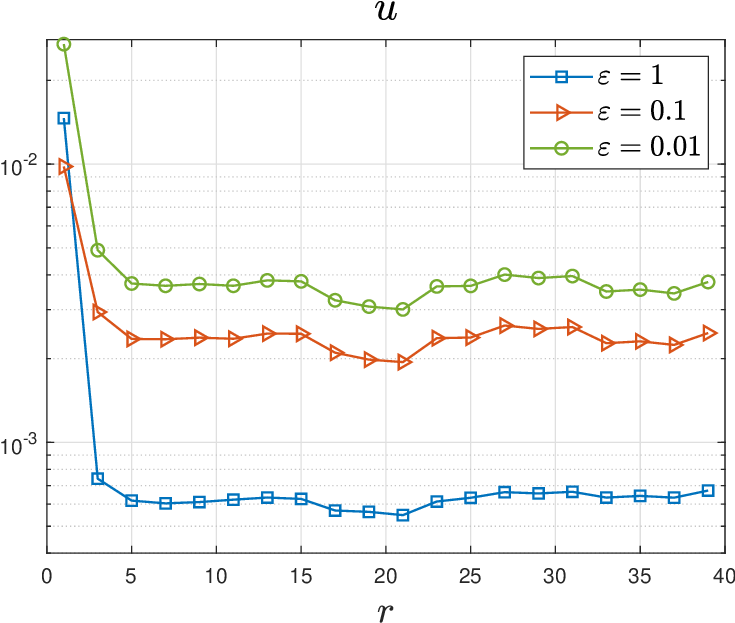}
	}       
	\caption{
		(UQ Test 2 in Subsection \ref{sec:uq_test2}) Errors of the bi-fidelity approximation for $n$ (left), $u$ (right) with respect to the number of high-fidelity simulation runs.}    
	\label{fig:uq2_error}
\end{figure}

\begin{figure}[!hptb]
	\centering
	\subfigure{
		\includegraphics[width = 0.4\textwidth]{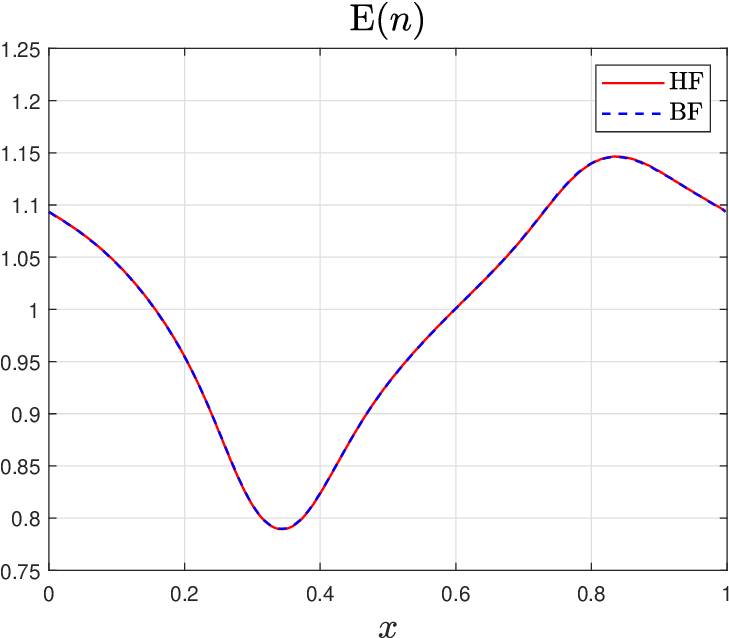}
	}
	\subfigure{               
		\includegraphics[width = 0.4\textwidth]{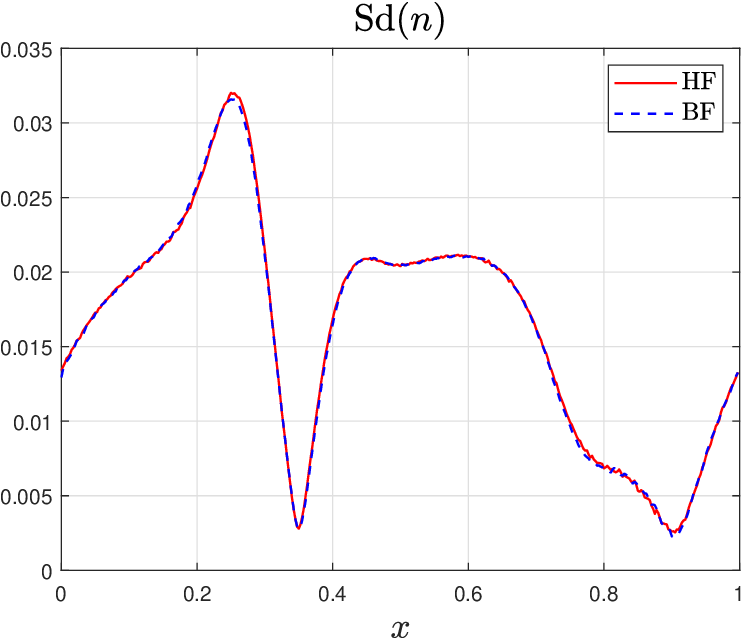}
	}       
	\subfigure{
		\includegraphics[width = 0.4\textwidth]{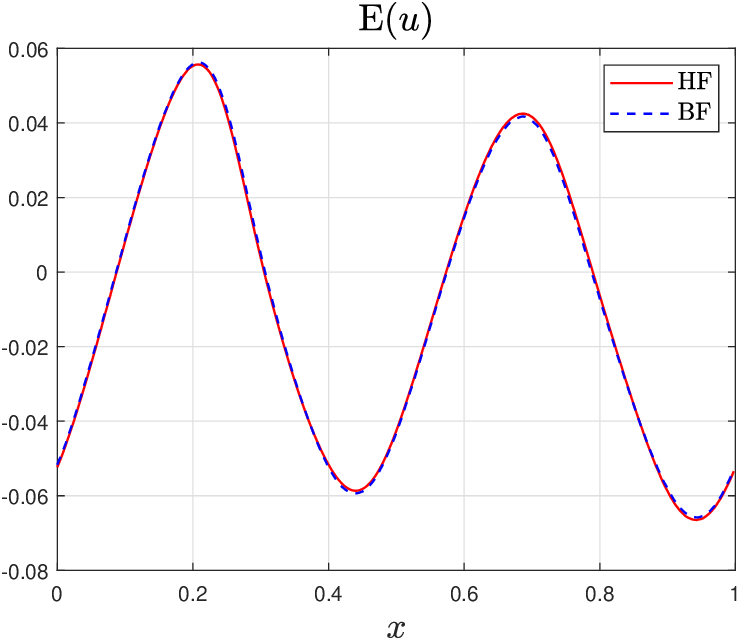}
	}
	\subfigure{               
		\includegraphics[width = 0.4\textwidth]{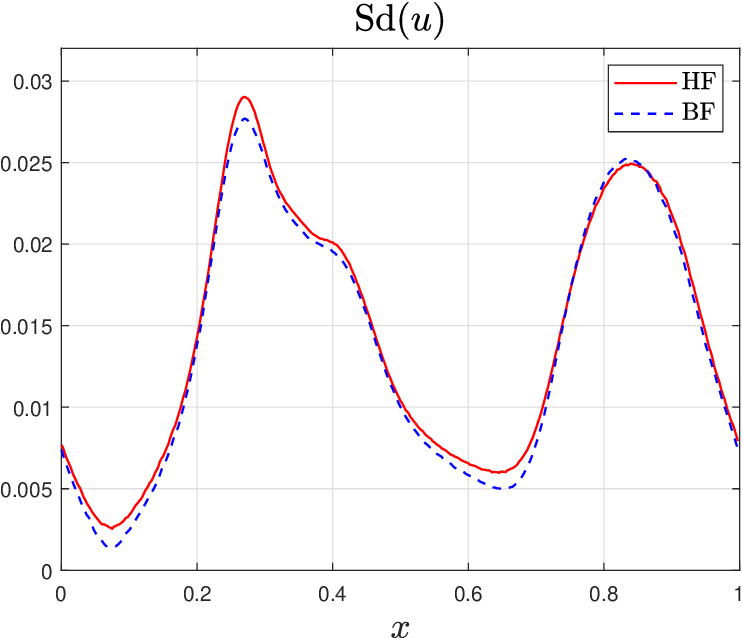}
	} 
	\caption{(UQ Test 2, $\e = 0.01$ in Subsection \ref{sec:uq_test2}) Mean and standard deviation of $n$, $u$ of high-fidelity and bi-fidelity solutions by using $r=10$.}    
	\label{fig:uq2_1e2_mean_std}
\end{figure}

In this test, we consider a more complicated initial distribution which contains uncertain parameters. Assume the initial condition given by 
\begin{equation}
	\label{eq:ini_uq2}
	f(0,x,v,\mathbf{z}) = n_{1}(x,\mathbf{z})\delta(v-u_{1}(x)) + n_{2}(x,\mathbf{z})\delta(v-u_{2}(x)),
\end{equation}
where $x\in\left[0,1\right]$, and
\begin{equation}
	\label{eq:ini_new_uq2}
	\begin{gathered}
		n_{1}(x,\mathbf{z})  = 1 - 0.6\cos^2{(2\pi x + \frac{\pi}{3})} + \sum_{i=1}^{5}\frac{1}{(i\pi)^2} \cos{(2\pi x)}z_i, \qquad 
		u_{1}(x) = 0.2\sin{(2\pi x)}, \\
		n_{2}(x,\mathbf{z})  = 1 + 0.4\sin{(2\pi x - \frac{\pi}{4})} + \sum_{i=1}^{5} \frac{1}{(i\pi)^2}\sin{(2\pi x)}z_i, \qquad 
		u_{2}(x) =  0.1\cos{(2\pi x)}. 
	\end{gathered}
\end{equation}
Here the random variable $\mathbf{z}= (z_1, \cdots, z_5)$ is 5-dimensional, with $\{z_i\}_{i=1}^5$ following the uniform distribution on $[-1,1]$. We set the final time $t=0.5$. In the HF model, $N_{h} = 400$, $N_{tol} = 1.6\times10^{6}$ and $N_c=600$ is used, with the time step the same as Test 1. The numerical parameters for the LF model are the same as used in Test 1. 
We notice a fast decay of $L^2$ errors between HF and BF solutions in Figure \ref{fig:uq2_error}, and the errors achieve a satisfactory level of accuracy around $O(10^{-3})$ when the number of HF runs is as small as $r=5$. 
This phenomenon stays uniformly true for problem in different regimes, with $\e$ ranging from 0.01, 0.1 to 1. 
In Figure \ref{fig:uq2_1e2_mean_std}, we test the model with $\e=0.01$ and show 
the mean and standard deviation of $n$ and $u$ by running the expensive HF solver only 10 times. One can observe that the BF approximations match really well with the HF solutions. Thus for our test 2 with more complex initial distribution, the same conclusion can be drawn as Test 1, that is, the BF model can effectively approximate the HF solution in different regimes, at a low computational cost. 

\subsubsection{Test 3: uncertain initial data with KL-expansion}
\label{sec:uq_test3}

In this last test, we consider a more realistic stochastic problem with uncertainties in the initial distribution. Let $f(0,x,v,\mathbf{z})$ be
\begin{equation}
	\label{eq:uq3_init_f}
	f(0,x,v,\mathbf{z}) = n_{0}(x,\mathbf{z})\delta(v-u_{0}(x)),\qquad 
	u_{0}(x) = \sin{(x)}.
\end{equation}
where $x\in \mathcal{X}$ and $\mathbf{z} \in \Omega$ represents a random sample in the probability space $\left(\Omega, \mathcal{F}, \mu\right)$. Here $n_{0}:\:\mathcal{X}\times \Omega \to \mathbb{R}$ is a square-integrable stochastic process, which can be decomposed in a
bi-orthogonal fashion through the Karhunen-Lo$\grave{e}$ve expansion \cite{bao2020inverse}. Assume 
\begin{equation}
	\label{eq:uq3_init_density}
	n_{0}(x,\mathbf{z}) = \tilde{n}(x) + \sum_{j=0}^{\infty}\sqrt{\lambda_{j}}z_{j}\varphi_{j}(x).
\end{equation}
with 
$\tilde{n}(x) = \frac{1}{s}\left(1.5 + 0.2\cos{(2x)}+ 0.1\cos{(4x)}\right)$, where $s$ is the normalization coefficient. The detailed expansion of \eqref{eq:uq3_init_density} is provided in Appendix \ref{sec:Annex-KL-expansion}.
We consider the initial distribution given by 
\begin{equation}
	\label{eq:uq3_init_f_n0}
	f(0,x,v,\mathbf{z}) = n_{0}(x,\mathbf{z})\delta(v-u_{0}(x)),\qquad u_{0}(x) = \sin{(x)}.
\end{equation}
Let the spatial domain $\mathcal{X}=[0,2\pi]$, the parameters in \eqref{eq:covar_operator_2} be $\sigma = \frac{1}{15}$, $l = 0.5$ with eigenvalues satisfying $\lambda_{j} > 10^{-6}$. For the HF model, $N_{h} = 400$, $ \Delta t_{h} = 0.4\frac{\Delta x}{\max\{u_0\}}$, $N_{tol} = 1.6\times10^{6}$ and $N_c=200$ repetition times is used to reduce the noise of PIC method.
In the LF solver, we take $N_{l} = 1000$ and $\Delta t_{l} = 0.1 \Delta x$. The model in different regimes ($\e = 1,0.1,0.01$) is studied, with final time $t = 0.5$. 

\begin{figure}[!hptb]
	\centering
	\subfigure{
		\includegraphics[width = 0.3\textwidth]{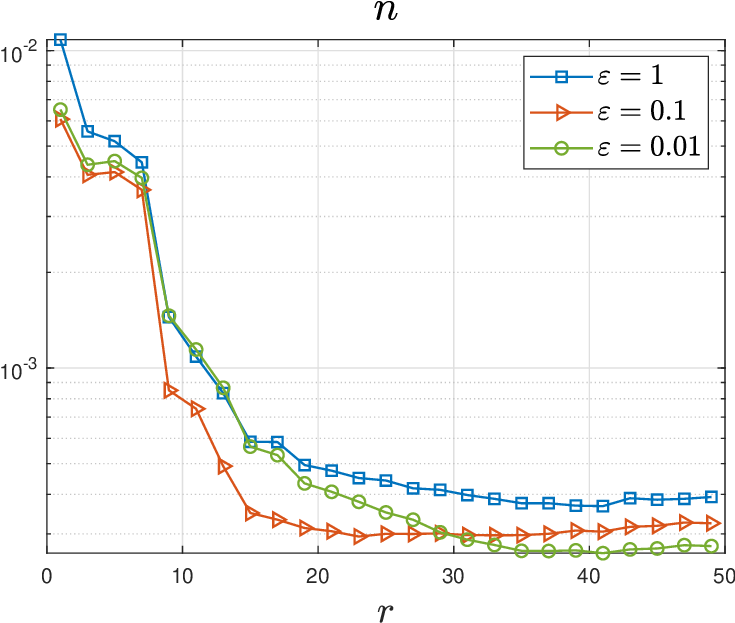}
	}
	\subfigure{               
		\includegraphics[width = 0.3\textwidth]{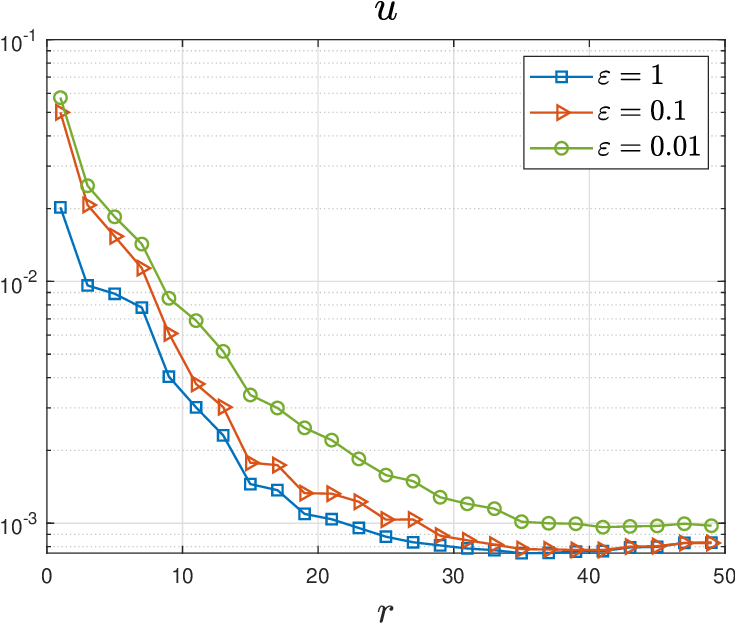}
	}       
	\caption{(UQ Test 3, $\e = 1,0.1,0.01$ in Figure \ref{sec:uq_test3}) Errors of the bi-fidelity approximation for $n$ (left), $u$ (right) with respect to the number of high-fidelity simulation runs.}  
	\label{fig:uq3_error}
\end{figure}

In Figure \ref{fig:uq3_error}, we detect a simlar convergence trend as in the previous two tests, however the $l_2$ errors of the BF approximation decay slightly slower. For the model in different regimes ($\e=1, 0.1, 0.01$), the errors seem to be saturated when $r$ reaches about 30, at the accuracy level of nearly $O(10^{-4})$ for the density $n$ and $O(10^{-3})$ for the velocity. In this test, we employ the KL expansion to characterize the initial distribution that contains uncertainties, which is more complicated than the settings in previous two tests. 
As a result, for the BF construction one needs more information to capture the HF solution space, that is, through utilizing larger number of HF implementations one is able to approximate the HF model at the same accuracy level as in the previous tests.

\begin{figure}[!hptb]
	\centering
	\subfigure{
		\includegraphics[width = 0.4\textwidth]{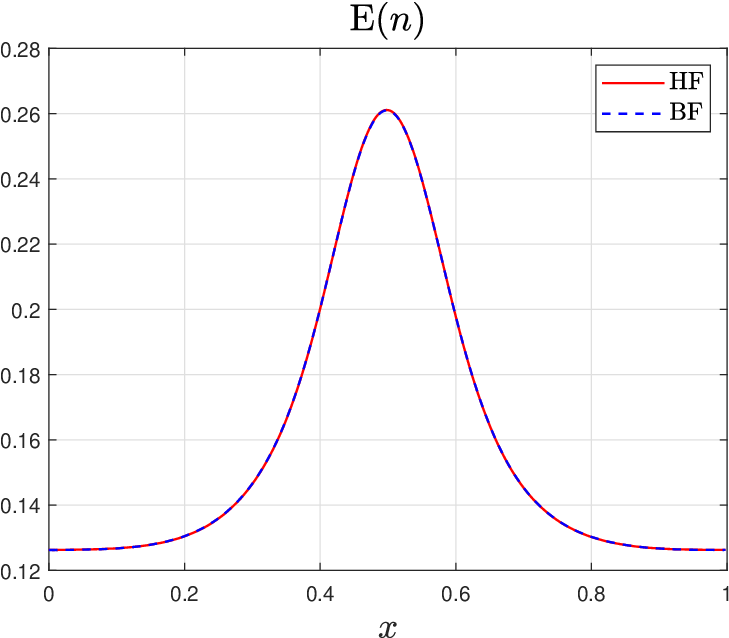}
	}
	\subfigure{               
		\includegraphics[width = 0.4\textwidth]{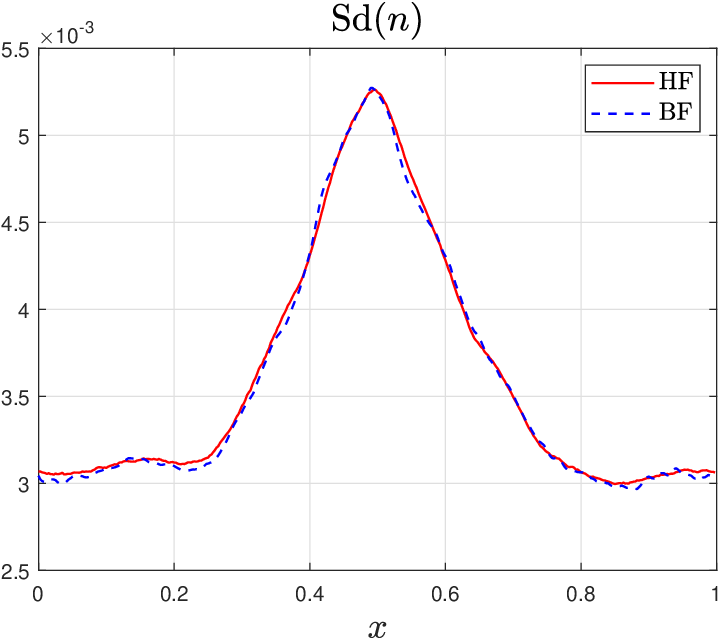}
	}       
	\subfigure{
		\includegraphics[width = 0.4\textwidth]{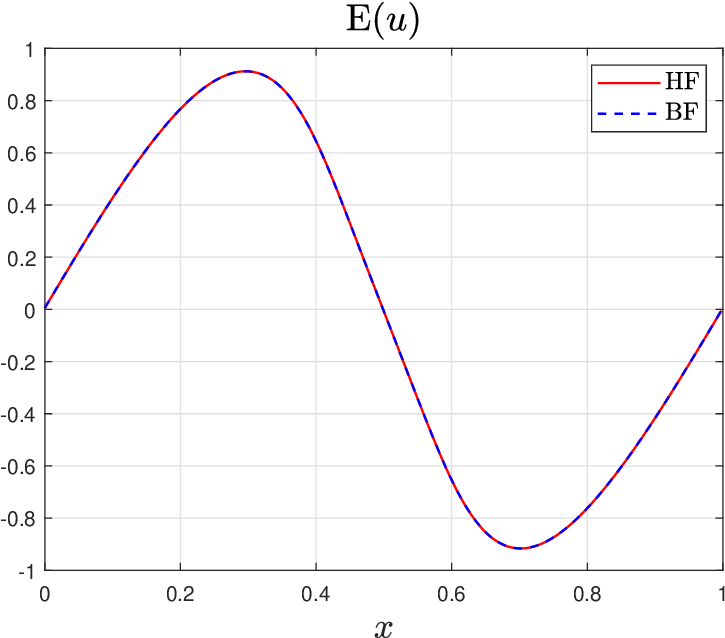}
	}
	\subfigure{               
		\includegraphics[width = 0.4\textwidth]{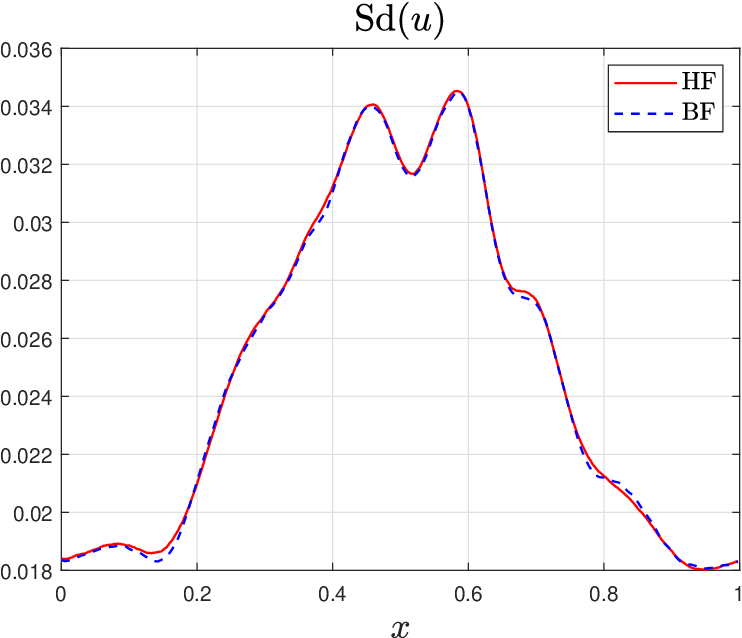}
	} 
	\caption{
		(UQ Test 3, $\e=0.1$ in Subsection \ref{sec:uq_test3}) Mean and standard deviation of $n$, $u$ of high-fidelity solutions and bi-fidelity solutions by using $r=30$.}    
	\label{fig:uq3_1e1_mean_std}
\end{figure}

Moreover, in Figure \ref{fig:uq3_1e1_mean_std} we present the expectation (left) and standard deviation (right) for quantities $n$, $u$ by using $r=30$ HF simulation runs in the BF approximation. It is clearly seen that the HF and BF solutions match quite well, thus demonstrating the efficiency and accuracy of our bi-fidelity method. 
\FloatBarrier 
\section{Conclusion}
\label{sec:conclusion}
	In this paper, inspired by the work \cite{degond2010asymptotic} we developed a AP-PIC method for solving the Vlasov-Poisson system with massless electrons where the Poisson equation is nonlinear, thus leading to new challenges. Our AP-PIC scheme is stable and accurate, it can capture the quasineutral limit system numerically, without resolving the space and time discretizations associated with the small Debye length in plasma. In the second part of the paper, we study the corresponding UQ problem and design an efficient bi-fidelity method in the stochastic collocation framework. By leveraging a low-fidelity model chosen by the macroscopic equations, our bi-fidelity approximations are capable of describing behaviour of the high-fidelity (VPME) model in the random space. A series of numerical experiments have validated the effectiveness and high efficiency of our proposed schemes. In the future work, we will tackle more complex kinetic models in plasmas, for example the problem with three-dimensional space and velocity variables, or model with higher-dimensional random uncertainties, in addition to design high-order AP schemes by incorporating the explicit time discretization mentioned in Remark \ref{remark: penalty}.
\section*{Acknowledgements}
\begin{acknowledgement}
L. Liu acknowledges the support by National Key R\&D Program of China (2021YFA1001200), Ministry of Science and Technology in China, Early Career Scheme (24301021) and General Research Fund (14303022 \& 14301423) funded by Research Grants Council of Hong Kong. This work of Yanli Wang is partially supported by the Foundation of the President of China Academy of Engineering Physics (YZJJZQ2022017), and the National Natural Science Foundation of China (Grant No. 12171026, U2230402, and 12031013). We also thank Dr. Yiwen Lin (Shanghai Jiao Tong University) on helpful discussions for the Karhunen-Lo\'eve expansions in our last numerical test. 
\end{acknowledgement}
\bibliographystyle{plain}
\bibliography{reference.bib}
\appendix
\section{Derivation of the model, classical PIC method and Karhunen-Lo\texorpdfstring{$\grave{e}$}{\grave{e}}ve Expansions}

\subsection{Derivation of isothermal compressible Euler system }

We gives the details of deriving the isothermal compressible Euler system given by equation (2.4) here. 
Taking the first and second moments in velocity for equation (2.1a), one gets the mass and momentum conservation equations given by
\begin{equation}
    \label{eq:moment}
    	\left\{\begin{array}{l}
		\partial_{t} n+\nabla_{\mathbf{x}} \cdot(n \mathbf{u})=0 ,\\[6pt]
		\partial_{t} (n\mathbf{u})+ \nabla_{\mathbf{x}} \cdot (n\mathbf{u} \otimes \mathbf{u})+ \nabla\!\cdot P = -\, n \nabla\phi. 
	\end{array}\right.
\end{equation}
Here $P$ is the non-dimensional ion pressure tensor. If we assume isotropy $P = pI$ (with $I$ being the identity matrix), then $\nabla\!\cdot P  = \nabla p$. The usual procedure to close this system is to assume a thermodynamic equation that describes the relation between $p$ and $n$, which is shown as 
\begin{equation}
    \label{eq:pressure}
    p= \tau n,\qquad\tau = \frac{T_i}{T_e}.
\end{equation}
This closure assumption for the ion moment system can be found in \cite[Section 2]{lieberman1994principles}. Based on this closed system, one considers the Euler-Poisson-Boltzmann (EPB) model \cite{P2012Numerical} by neglecting the ion pressure tensor, where the cold-ion approximation (i.e., $T_i \to 0$, so $\tau \ll 1$) is adopted: 
\begin{equation}
    \label{eq:moment_nopress}
    	\left\{\begin{array}{l}
		\partial_{t} n+\nabla_{\mathbf{x}} \cdot(n \mathbf{u})=0 ,\\[6pt]
		\partial_{t} (n\mathbf{u})+ \nabla_{\mathbf{x}} \cdot (n\mathbf{u} \otimes \mathbf{u}) = -\, n \nabla\phi. 
	\end{array}\right.
\end{equation}
If the quasi-neutral assumption is made, the Poisson equation (2.1b) is replaced by the constraint of local charge being zero: 
\begin{equation*}
    n = e^{\phi}. 
\end{equation*}
Then one can write
\begin{equation*}
    n\nabla\phi = e^{\phi}\nabla\phi = \nabla e^{\phi} = \nabla n,
\end{equation*}
and the quasi-neutral Euler-Poisson-Boltzmann model is exactly the Isothermal Compressible Euler model given in equation (2.4): 
\begin{equation}
    	\left\{\begin{array}{l}
		\partial_{t} n+\nabla_{\mathbf{x}} \cdot(n \mathbf{u})=0 ,\\[6pt]
		\partial_{t} (n\mathbf{u})+ \nabla_{\mathbf{x}} \cdot (n\mathbf{u} \otimes \mathbf{u})= -\nabla n. 
	\end{array}\right.
\end{equation}

\subsection{Details of the classical PIC method}
Here we mention the assignment of weights of particles and some details in the classical PIC method.
The Nearest Grid Point (NGP) model is adopted to determine the initial charge density assignment $w_{k}$ in \eqref{eq:f-PIC}.
\begin{equation}
	\label{eq:pic-charge-assignment}
	\begin{aligned}
		n^{0}_{j} &  = \frac{1}{h}\int_{x_{j}-\frac{h}{2}}^{x_{j}+\frac{h}{2}}\mathrm{d}x\:\int_{-\infty}^{+\infty}f^{0}(x,v)\mathrm{d}v \\
		& \approx \frac{1}{h}\int_{x_{j}-\frac{h}{2}}^{x_{j}+\frac{h}{2}}\mathrm{d}x\:\int_{-\infty}^{+\infty}f^{0}_{N_{\text{total}}}(x,v)\mathrm{d}v.
	\end{aligned}
\end{equation}
where $n^{0}_{j}$ is the initial average charge density in the volume $\left[x_{j}-\frac{h}{2}, x_{j}+\frac{h}{2}\right]$.
In the NGP model, the charge assignment function $W$ is introduced to deal with the term $\delta(x-X_{k})$ which is given by
\begin{equation}
	\label{Annex:pic-assignment-function}
	W(x)=\left\{\begin{array}{ll}
		1 & |x|<h / 2 \text { or } x=h / 2, \\[4pt]
		0 & \text { otherwise }
	\end{array}\right.
\end{equation}
Then the weight $w_{k}$ can be derived as follows 
\begin{equation}
	\label{eq:pic-weight-1}
	\frac{1}{h}\sum_{\substack{\text { particles k} \\ \text { in cell } j}}w_{k} = n^{0}_{j}.
\end{equation}
For convenience, the initial number of particles per cell is set as the same, which depends on the total number of particles $N_{\text{total}}$ and the size of the mesh $N_{h}$. Besides, the weights of particles in same cell are set as the same initially, which means for all particles in cell $j$ at $t=0$, their weights can be derived by \eqref{eq:pic-weight-1}
\begin{equation}
	\label{eq:pic-weight-2}
	w_{k} = \frac{hn^{0}_{j}}{M},\qquad M = \frac{N_{\text{total}}}{N_{h}}.
\end{equation}
where $k$ represent particles in cell $j$. For detailed proof and analysis of the validity of the assignment function $W$, we refer to Chapter 5 of \cite{hockney2021computer}. 
\subsection{Karhunen-Lo\texorpdfstring{$\grave{e}$}{\grave{e}}ve Expansions}
\label{sec:Annex-KL-expansion}
Here we refer to some text books \cite{lord2014introduction, sullivan2015introduction} for details on the Karhunen-Lo$\grave{e}$ve expansion. We mention the following theorem: 
\begin{thrm}
	\label{thm:KL-expansion}
	(Karhunen-Lo$\grave{e}$ve). Let $U:\: \mathcal{X} \times \Omega \to \mathbb{R}$ be square-integrable stochastic process, with mean zero and continuous and square-integrable covariance function. Then
	\begin{equation}
		\label{eq:KL-expansion-thm-1}
		U = \sum_{j\in\mathbb{N}}Z_{j}\varphi_{j}.
	\end{equation}
	in $L^2$, where the $\{\varphi_j\}_{j\in\mathbb{N}}$ are orthonormal eigenfunctions of the covariance operator $C_{U}$, the corresponding eigenvalues $\{\lambda_j\}_{j\in\mathbb{N}} $ are non-negative, the convergence of the series is in $L^{2}(\Omega, \mu; \mathbb{R})$ and uniform among compact families of $x\in\mathcal{X}$, with
	\begin{equation}
		\label{eq:KL-expansion-thm-2}
		Z_{j} = \int_{\mathcal{X}}U(x)\varphi_{j}(x)\mathrm{d}x.
	\end{equation}
	Furthermore, the random variables $Z_{j}$ are centered, uncorrelated, and have variance $\lambda_{j}$:
	\begin{equation}
		\label{eq:KL-expansion-thm-3}
		\mathbb{E}_{\mu}\left[Z_{j}\right] = 0,\qquad \mathbb{E}_{\mu}\left[Z_{j}Z_{k}\right] = \lambda_{j}\delta_{jk}.
	\end{equation}
\end{thrm}
The above theorem demonstrates the validity of the Karhunen-Lo$\grave{e}$ve expansion and how to compute the eigenvalues after determining the orthogonal basis. Details of the proof and properties of $U$ can be found in Chapter 11 of \cite{sullivan2015introduction}. Regarding the initial density \eqref{eq:uq3_init_density} in Subsection \ref{sec:uq_test3}, which is equivalent to $U$ in the above theorem, we denote $x\in [0,\Lambda]$ with $\Lambda=2\pi$ and compute the covariance function \cite{bao2020inverse}:
\begin{equation}
	\label{eq:covar_operator_2}
	c(x-y) = \sigma^2 e^{-\frac{\left|x-y\right|^2}{l^{2}}},\: 0<l\ll \Lambda.
\end{equation}
where $\sigma$ is the root mean square of the stochastic process and $l$ is the correlation length. The orthonormal eigenfunctions of the covariance operator are defined by
\begin{equation}
	\label{eq:KL-ex-basis}
	\varphi_{j}(x)=\left\{\begin{array}{ll}
		\sqrt{\frac{1}{\Lambda}}, & j=0 \\
		\sqrt{\frac{2}{\Lambda}} \cos \left(\frac{2 j \pi x}{\Lambda}\right), & j>1, \text { even } \\
		\sqrt{\frac{2}{\Lambda}} \sin \left(\frac{2 j \pi x}{\Lambda}\right), & j>1, \text { odd }
	\end{array}\right.
\end{equation}
The corresponding eigenvalues $\{\lambda_{j}\}_{j=0}^{\infty}$ are arranged in a descending order, which is computed by
\begin{equation}
	\label{eq:eigenvalue_compute}
	\lambda_{j} = \int_{-\pi}^{\pi}c(x)\varphi_{j}(x)\mathrm{d}x.
\end{equation}
and $\{z_{j}\}_{j=0}^{\infty},\: z_{j}\sim \mathcal{N}(0, 1)$. Note that the orthonormal basis \eqref{eq:KL-ex-basis} is not unique. Here, we choose the trigonometric function system.
\end{document}